\documentclass[12pt]{article}
\usepackage{amsmath,amssymb,amsthm,amscd,a4wide}

\begin{document}

\newcommand{\End}{{\rm{End}\ts}}
\newcommand{\Hom}{{\rm{Hom}}}
\newcommand{\ch}{{\rm{ch}\ts}}
\newcommand{\non}{\nonumber}
\newcommand{\wt}{\widetilde}
\newcommand{\wh}{\widehat}
\newcommand{\ot}{\otimes}
\newcommand{\la}{\lambda}
\newcommand{\La}{\Lambda}
\newcommand{\De}{\Delta}
\newcommand{\al}{\alpha}
\newcommand{\be}{\beta}
\newcommand{\ga}{\gamma}
\newcommand{\ka}{\kappa}
\newcommand{\si}{\sigma}
\newcommand{\vp}{\varphi}
\newcommand{\de}{\delta^{}}
\newcommand{\ze}{\zeta}
\newcommand{\om}{\omega}
\newcommand{\hra}{\hookrightarrow}
\newcommand{\ve}{\varepsilon}
\newcommand{\ts}{\,}
\newcommand{\vac}{\mathbf{1}}
\newcommand{\di}{\partial}
\newcommand{\qin}{q^{-1}}
\newcommand{\tss}{\hspace{1pt}}
\newcommand{\Sr}{ {\rm S}}
\newcommand{\U}{ {\rm U}}
\newcommand{\Y}{ {\rm Y}}
\newcommand{\AAb}{\mathbb{A}\tss}
\newcommand{\CC}{\mathbb{C}\tss}
\newcommand{\QQ}{\mathbb{Q}\tss}
\newcommand{\SSb}{\mathbb{S}\tss}
\newcommand{\ZZ}{\mathbb{Z}\tss}
\newcommand{\Z}{{\rm Z}}
\newcommand{\Ac}{\mathcal{A}}
\newcommand{\Lc}{\mathcal{L}}
\newcommand{\Mc}{\mathcal{M}}
\newcommand{\Pc}{\mathcal{P}}
\newcommand{\Qc}{\mathcal{Q}}
\newcommand{\Tc}{\mathcal{T}}
\newcommand{\Sc}{\mathcal{S}}
\newcommand{\Bc}{\mathcal{B}}
\newcommand{\Ec}{\mathcal{E}}
\newcommand{\Hc}{\mathcal{H}}
\newcommand{\Uc}{\mathcal{U}}
\newcommand{\Wc}{\mathcal{W}}
\newcommand{\Ar}{{\rm A}}
\newcommand{\Ir}{{\rm I}}
\newcommand{\Zr}{{\rm Z}}
\newcommand{\gl}{\mathfrak{gl}}
\newcommand{\Pf}{{\rm Pf}}
\newcommand{\oa}{\mathfrak{o}}
\newcommand{\spa}{\mathfrak{sp}}
\newcommand{\g}{\mathfrak{g}}
\newcommand{\h}{\mathfrak h}
\newcommand{\n}{\mathfrak n}
\newcommand{\z}{\mathfrak{z}}
\newcommand{\Zgot}{\mathfrak{Z}}
\newcommand{\p}{\mathfrak{p}}
\newcommand{\sll}{\mathfrak{sl}}
\newcommand{\agot}{\mathfrak{a}}
\newcommand{\qdet}{ {\rm qdet}\ts}
\newcommand{\Ber}{ {\rm Ber}\ts}
\newcommand{\HC}{ {\mathcal HC}}
\newcommand{\cdet}{ {\rm cdet}}
\newcommand{\tr}{ {\rm tr}}
\newcommand{\str}{ {\rm str}}
\newcommand{\loc}{{\rm loc}}
\newcommand{\Gr}{ {\rm Gr}\tss}
\newcommand{\sgn}{ {\rm sgn}\ts}
\newcommand{\ba}{\bar{a}}
\newcommand{\bb}{\bar{b}}
\newcommand{\bi}{\bar{\imath}}
\newcommand{\bj}{\bar{\jmath}}
\newcommand{\bk}{\bar{k}}
\newcommand{\bl}{\bar{l}}
\newcommand{\Sym}{\mathfrak S}
\newcommand{\fand}{\quad\text{and}\quad}
\newcommand{\Fand}{\qquad\text{and}\qquad}
\newcommand{\vt}{{\tss|\hspace{-1.5pt}|\tss}}

\renewcommand{\theequation}{\arabic{section}.\arabic{equation}}

\newtheorem{thm}{Theorem}[section]
\newtheorem{lem}[thm]{Lemma}
\newtheorem{prop}[thm]{Proposition}
\newtheorem{cor}[thm]{Corollary}
\newtheorem{conj}[thm]{Conjecture}

\theoremstyle{definition}
\newtheorem{defin}[thm]{Definition}

\theoremstyle{remark}
\newtheorem{remark}[thm]{Remark}
\newtheorem{example}[thm]{Example}

\newcommand{\bth}{\begin{thm}}
\renewcommand{\eth}{\end{thm}}
\newcommand{\bpr}{\begin{prop}}
\newcommand{\epr}{\end{prop}}
\newcommand{\ble}{\begin{lem}}
\newcommand{\ele}{\end{lem}}
\newcommand{\bco}{\begin{cor}}
\newcommand{\eco}{\end{cor}}
\newcommand{\bde}{\begin{defin}}
\newcommand{\ede}{\end{defin}}
\newcommand{\bex}{\begin{example}}
\newcommand{\eex}{\end{example}}
\newcommand{\bre}{\begin{remark}}
\newcommand{\ere}{\end{remark}}
\newcommand{\bcj}{\begin{conj}}
\newcommand{\ecj}{\end{conj}}

\newcommand{\bal}{\begin{aligned}}
\newcommand{\eal}{\end{aligned}}
\newcommand{\beq}{\begin{equation}}
\newcommand{\eeq}{\end{equation}}
\newcommand{\ben}{\begin{equation*}}
\newcommand{\een}{\end{equation*}}

\newcommand{\bpf}{\begin{proof}}
\newcommand{\epf}{\end{proof}}

\def\beql#1{\begin{equation}\label{#1}}

\title{\Large\bf The MacMahon Master Theorem for right quantum
superalgebras and higher Sugawara operators for $\wh\gl_{m|n}$}

\author{{A. I. Molev\quad and\quad E. Ragoucy}}

\date{} 
\maketitle

\vspace{5 mm}

\begin{abstract}
We prove an analogue of the MacMahon Master Theorem for the right
quantum superalgebras. In particular, we obtain a new and simple
proof of this theorem for the right quantum algebras. In the super
case the theorem is then used to construct higher order Sugawara
operators for the affine Lie superalgebra $\wh\gl_{m|n}$ in an
explicit form. The operators are elements of a completed universal
enveloping algebra of $\wh\gl_{m|n}$ at the critical level. They
occur as the coefficients in the expansion of a noncommutative
Berezinian and as the traces of powers of generator matrices.
The same construction yields higher Hamiltonians for the Gaudin model
associated with the Lie superalgebra $\gl_{m|n}$.
We also use the Sugawara operators to produce
algebraically independent generators of the algebra of singular
vectors of any generic Verma module at the critical level
over the affine Lie superalgebra.
%
\end{abstract}

\vspace{5 mm}

\vspace{25 mm}

\noindent
School of Mathematics and Statistics\newline
University of Sydney,
NSW 2006, Australia\newline
alexm@maths.usyd.edu.au

\vspace{7 mm}

\noindent
LAPTH, Chemin de Bellevue, BP 110\newline
F-74941 Annecy-le-Vieux cedex, France\newline
ragoucy@lapp.in2p3.fr

\newpage

\hfill {\sl Dedicated to Grigori Olshanski on his 60th birthday}

\section{Introduction}\label{sec:int}
\setcounter{equation}{0}

{\bf 1.1. MacMahon Master Theorem.}\quad A natural
quantum analogue of the celebrated
MacMahon Master Theorem was proved
by Garoufalidis, L\^{e} and Zeilberger in \cite{glz:qm}.
A few different proofs and generalizations
of this analogue have since been
found; see Foata and Han~\cite{fh:np,fh:br},
Hai and Lorenz~\cite{hl:ka} and
Konvalinka and Pak~\cite{kp:ne}.
In the quantum MacMahon Master Theorem of \cite{glz:qm}
the numerical matrices are replaced with the {\it right quantum
matrices\/} $Z=[z_{ij}]$ whose matrix elements
satisfy some quadratic relations involving a parameter $q$.
As explained in \cite{fh:br}, this parameter may be
taken to be equal to $1$ without a real loss of generality.
Then the relations for the matrix elements
$z_{ij}$ take the form
\beql{rq}
[z_{ij}, z_{kl}]=[z_{kj}, z_{il}]
\qquad\text{for all}\quad i,j,k,l\in \{1,\dots,N\},
\eeq
where $[x,y]=xy-yx$.

The proof of the quantum MacMahon Master Theorem given in
\cite{hl:ka} is based on the theory of Koszul algebras.
In that approach, \eqref{rq} are the defining relations
of the bialgebra $\underline{\rm end}\ts \Ac$ associated
with the symmetric algebra $\Ac$ of a vector space.
The construction of the bialgebra $\underline{\rm end}\ts \Ac$
associated with an arbitrary quadratic algebra
or superalgebra $\Ac$
is due to Manin~\cite{m:sr, m:qg}, and the matrices $Z$
satisfying \eqref{rq} are also known as {\it Manin matrices\/};
see \cite{css:nd}, \cite{cf:mm}, \cite{cfr:ap}.
In the super case the construction applied to
the symmetric algebra of a $\ZZ_2$-graded vector
space (i.e., superspace) of dimension $m+n$
leads to the defining relations of the {\it right quantum
superalgebra\/} $\Mc_{m|n}$. This superalgebra
is generated by elements $z_{ij}$ with
$\ZZ_2$-degree (or parity)
$\bi+\bj$, where $\bi=0$ for $1\leqslant i\leqslant m$ and
$\bi=1$ for $m+1\leqslant i\leqslant m+n$.
The defining relations have the form
\beql{srq}
[z_{ij}, z_{kl}]=[z_{kj}, z_{il}](-1)^{\bi\bj+\bi\bk+\bj\bk}
\qquad\text{for all}\quad i,j,k,l\in \{1,\dots,m+n\},
\eeq
where $[x,y]=xy-yx(-1)^{\deg x\deg y}$ is the super-commutator
of homogeneous elements $x$ and $y$,
as presented e.g. in
\cite[Example~3.14]{hkl:nh}\footnote{The relations
given in \cite{hkl:nh} correspond to our
{\it left\/} quantum superalgebra as defined
in Sec.~\ref{subsec:arqs} below.}.
We will call any matrix $Z=[z_{ij}]$ satisfying
\eqref{srq} a Manin matrix.

Our first main result (Theorem~\ref{thm:mmmt})
is an analogue of the MacMahon Master Theorem
for the right quantum superalgebra $\Mc_{m|n}$.
We will identify the matrix $Z$ with an element
of the tensor product superalgebra $\End\CC^{m|n}\ot \Mc_{m|n}$ by
\ben
Z=\sum_{i,j=1}^{m+n}e_{ij}\ot z_{ij}(-1)^{\bi\tss\bj+\bj},
\een
where the $e_{ij}$ denote the standard matrix units.
Taking multiple tensor products
\beql{tenpr}
\End\CC^{m|n}\ot\dots\ot\End\CC^{m|n}\ot\Mc_{m|n}
\eeq
with $k$ copies of $\End\CC^{m|n}$,
for any $a=1,\dots,k$ we will write $Z_a$
for the matrix $Z$ corresponding to the $a$-th copy of
the endomorphism superalgebra
so that the components in all remaining
copies are the identity matrices. The symmetric group $\Sym_k$
acts naturally on the tensor product
space $(\CC^{m|n})^{\ot k}$. We let $H_k$ and $A_k$ denote
the respective images of the normalized
symmetrizer and antisymmetrizer
\beql{symm}
\frac{1}{k!}\sum_{\si\in\Sym_k}\si\in\CC[\Sym_k],
\qquad\qquad
\frac{1}{k!}\sum_{\si\in\Sym_k}\sgn\si\cdot\si\in\CC[\Sym_k]
\eeq
in \eqref{tenpr}. Recall that the {\it supertrace\/}
of an even matrix $X=[x_{ij}]$ is defined by
\beql{supertr}
\str\ts X
=\sum_{i=1}^{m+n} x_{ii}(-1)^{\bi}.
\eeq
Now set
\beql{bf}
{\rm Bos}=1+\sum_{k=1}^{\infty} \str\ts H_k Z_1\dots Z_k,
\qquad
{\rm Ferm}=1+\sum_{k=1}^{\infty} (-1)^k\ts\str\ts A_k Z_1\dots Z_k,
\eeq
where $\str$ denotes the supertrace taken with respect to
all $k$ copies of $\End\CC^{m|n}$. Our analogue of the
MacMahon Master Theorem for the right quantum superalgebra
$\Mc_{m|n}$ now reads
\beql{mmmt}
{\rm Bos}\times {\rm Ferm}=1.
\eeq

The summands in \eqref{bf} can be written explicitly
in terms of the generators $z_{ij}$;
see Proposition~\ref{prop:extsh} below.
In each of the particular cases $n=0$ and $m=0$ the identity
\eqref{mmmt}
turns into the quantum MacMahon Master Theorem of \cite{glz:qm}.

Our proof of the identity \eqref{mmmt} is based on the use
of the matrix form of the defining relations \eqref{srq}.
These relations can be written as
\beql{mdefrel}
(1-P_{12})\ts [Z_1,Z_2]=0
\eeq
which is considered as an identity in the superalgebra
\eqref{tenpr} with $k=2$ and $P_{12}$ is the image
of the transposition $(12)\in\Sym_2$. The proof of \eqref{mmmt}
is derived from \eqref{mdefrel} by using some elementary properties
of the symmetrizers and antisymmetrizers.

\medskip
\noindent
{\bf 1.2. Noncommutative Berezinian.}\quad It is well-known that
in the super-commutative specialization
the expressions \eqref{bf}
coincide with
the expansions of the {\it Berezinian\/} so that
\ben
{\rm Ferm}=\Ber(1-Z)\Fand {\rm Bos}=\big[\Ber(1-Z)\big]^{-1};
\een
see \cite{b:is}, \cite{kv:be}, \cite{s:si}.
In order to prove their noncommutative analogues,
we introduce the {\it affine
right quantum superalgebra\/} $\wh\Mc_{m|n}$.
It is generated by a countable set of elements $z^{(r)}_{ij}$
of parity $\bi+\bj$, where $r$ runs over the set of positive
integers. The defining relations of $\wh\Mc_{m|n}$
take the form
\beql{mdefrelaff}
(1-P_{12})\ts [Z_1(u),Z_2(u)]=0,
\eeq
where $u$ is a formal variable, and the matrix elements
of the matrix $Z(u)=[z_{ij}(u)]$ are the formal power series
\beql{serz}
z_{ij}(u)=\de_{ij}+z^{(1)}_{ij}\ts u+z^{(2)}_{ij}\ts u^2+\cdots.
\eeq
The relation \eqref{mdefrelaff} is understood in the sense that
all coefficients of the powers of $u$ on the left hand side vanish.
The formal power series $Z(u)$ in $u$ with matrix coefficients
is invertible and
we denote by $z^{\tss\prime}_{ij}(u)$ the matrix elements of its inverse
so that $Z^{-1}(u)=[z^{\tss\prime}_{ij}(u)]$. Now we set
\begin{align}
\Ber Z(u){}&=\sum_{\si\in\Sym_m}\sgn\si\cdot
z_{\si(1)1}(u)\dots z_{\si(m)m}(u)
\non\\
\label{ber}
{}&\times{}
\sum_{\tau\in\Sym_n}\sgn\tau\cdot
z^{\tss\prime}_{m+1,m+\tau(1)}(u)\dots z^{\tss\prime}_{m+n,m+\tau(n)}(u).
\end{align}
In the super-commutative specialization $\Ber Z(u)$
coincides with the ordinary Berezinian of the matrix $Z(u)$;
see \cite{b:is}.

We prove the following expansions of the noncommutative Berezinian,
where $Z$ is a Manin matrix:
\beql{charferm}
\Ber(1+uZ)=\sum_{k=0}^{\infty} u^k\ts\str\ts A_k Z_1\dots Z_k
\eeq
and hence by \eqref{mmmt}
\beql{charbos}
\big[\Ber(1-uZ)\big]^{-1}=\sum_{k=0}^{\infty}
u^k\ts\str\ts H_k Z_1\dots Z_k.
\eeq
Furthermore, we derive the following super-analogues
of the Newton identities:
\beql{newton}
\big[\Ber(1+uZ)\big]^{-1}\ts\di_u\tss \Ber(1+uZ)
=\sum_{k=0}^{\infty} (-u)^k\ts\str\ts Z^{k+1}
\eeq
which allow one to express the elements $\str\ts Z^{\tss k}$
in terms of the coefficients of either series \eqref{charferm}
or \eqref{charbos}.
In the even case of $n=0$ the Newton identities were proved
in the papers \cite{cf:mm, cfr:ap} which also contain
various generalizations of the matrix algebra
properties to the class of Manin matrices.
However, the proof of the identities
given in these papers relies on the existence of the
adjoint matrix and does not immediately
extend to the super case. To prove \eqref{newton}
we employ instead an appropriate
super-extension of the $R$-matrix arguments used in \cite{iop:gc};
see also \cite{gps:sp}.

\medskip
\noindent
{\bf 1.3. Segal--Sugawara vectors.}\quad The above properties of
Manin matrices will be used in our construction
of the {\it higher order Sugawara operators}
for the affine Lie superalgebra $\wh\gl_{m|n}$.
The commutation relations of the Lie superalgebra
$\wh\gl_{m|n}=\gl_{m|n}[t,t^{-1}]\oplus\CC K$
have the form
\begin{align}
\non
\big[e_{ij}[r],e_{kl}[s\tss]\tss\big]
=\de_{kj}\ts e_{i\tss l}[r+s\tss]
{}&-\de_{i\tss l}\ts e_{kj}[r+s\tss](-1)^{(\bi+\bj)(\bk+\bl)}\\
\label{commrel}
{}&+K\Big(\de_{kj}\tss\de_{i\tss l}(-1)^{\bi}
-\frac{\de_{ij}\tss\de_{kl}}{m-n}(-1)^{\bi+\bk}\Big)
\ts r\tss\de_{r,-s}
\end{align}
(assuming $m\ne n$), the element $K$ is even and central, and
we set $e_{ij}[r]=e_{ij}t^r$. To include the case
$m=n$ we will re-scale the central element by setting $K=K'(n-m)$
so that \eqref{commrel} will be written uniformly for all values of $m$ and $n$.
We will also consider
the extended Lie superalgebra $\wh\gl_{m|n}\oplus\CC\tau$,
where the element $\tau$ is even and
\beql{tauco}
\big[\tau,e_{ij}[r]\tss\big]=-r\ts e_{ij}[r-1],\qquad
\big[\tau,K\big]=0.
\eeq

For any $\ka\in\CC$ the {\it affine vertex algebra}
$V_{\ka}(\gl_{m|n})$ can be defined as the quotient
of the universal enveloping algebra
$\U(\wh\gl_{m|n})$ by the left ideal generated by
$\gl_{m|n}[t]$ and $K-\ka$; see e.g.
\cite{fb:va}, \cite{k:va}. The {\it center\/} of the vertex
algebra $V_{\ka}(\gl_{m|n})$ is its subspace
spanned by all elements $b\in V_{\ka}(\gl_{m|n})$
such that $\gl_{m|n}[t]\ts b=0$.
The axioms of the vertex algebra imply that the center
is a commutative associative superalgebra.
The center of the vertex algebra $V_{\ka}(\gl_{m|n})$ is trivial for
all values of $\ka$ except for the
{\it critical value\/} $\ka=n-m$ (this corresponds
to the evaluation $K'=1$ in terms of the re-scaled
central element $K'$).
In the latter
case the center is big and we denote it by
$\z(\wh\gl_{m|n})$. Any element of $\z(\wh\gl_{m|n})$
is called a {\it Segal--Sugawara vector\/}.
As a vector superspace, the vertex algebra $V_{n-m}(\gl_{m|n})$
can be identified with the universal enveloping algebra
$\U(t^{-1}\gl_{m|n}[t^{-1}])$. Moreover,
the multiplication in the commutative algebra
$\z(\wh\gl_{m|n})$ coincides with the multiplication in
$\U(t^{-1}\gl_{m|n}[t^{-1}])$. Therefore, the center
$\z(\wh\gl_{m|n})$ is naturally identified with
a commutative subalgebra of $\U(t^{-1}\gl_{m|n}[t^{-1}])$.

Another principal result of this paper
(Theorem~\ref{thm:traces} and Corollary~\ref{cor:ss})
is an explicit construction
of several families of Segal--Sugawara vectors.
The construction is based
on the observation that the matrix elements of the
matrix
\ben
\tau+\wh E[-1]=\big[\de_{ij}\tau+e_{ij}[-1](-1)^{\bi}\big]
\een
satisfy the defining relations \eqref{srq} of the right
quantum superalgebra. That is, $\tau+\wh E[-1]$ is a Manin
matrix. We use this fact to show that
all coefficients $s_{kl}$ in the expansion
\ben
\str(\tau+\wh E[-1])^k=
s_{k\tss 0}\ts\tau^{k}+s_{k1}\ts\tau^{k-1}
+\dots+s_{kk}
\een
are Segal--Sugawara vectors.
The relations
\eqref{charferm}, \eqref{charbos} and \eqref{newton}
then lead to explicit formulas for
a few other families of Segal--Sugawara vectors.
In particular, such vectors also occur
as the coefficients $b_{kl}$ in the expansion of
the Berezinian
\beql{expanber}
\Ber\big(1+u(\tau+\wh E[-1])\big)=
\sum_{k=0}^{\infty} \sum_{l=0}^k b_{kl}\ts u^k\tau^{k-l}.
\eeq

\medskip
\noindent
{\bf 1.4. Commutative subalgebras and higher Gaudin
Hamiltonians.}\quad Applying the state-field
correspondence map to the Segal--Sugawara vectors
of any of the families above,
we get explicit formulas
for elements of the center of the local completion
of the universal enveloping algebra $\U(\wh\gl_{m|n})$
at the critical level.
Such elements are called the (higher order) {\it Sugawara
operators\/}.
In particular, such operators can be calculated as
the Fourier coefficients of the fields $b_{kl}(z)$
defined by the expansion of
the normally ordered Berezinian
\beql{berfurier}
:\Ber\big(1+u(\di_z+\wh E(z))\big):{}=
\sum_{k=0}^{\infty} \sum_{l=0}^k b_{kl}(z)\ts u^k\di_z^{\ts k-l},
\eeq
where $\wh E(z)=[e_{ij}(z)(-1)^{\bi}]$ and
\ben
e_{ij}(z)=\sum_{r\in\ZZ} e_{ij}[r]\ts z^{-r-1}.
\een
Similarly, all Fourier coefficients of the fields
$s_{kl}(z)$ defined by the expansion
of the normally ordered supertrace
\beql{expanstr}
:\str\big(\di_z+\wh E(z)\big)^k:{}
= s_{k\tss 0}(z)\ts\di_z^{\tss k}+s_{k1}(z)\ts\di_z^{\tss k-1} +\dots
+s_{kk}(z)
\eeq
are Sugawara operators.

It was observed in \cite{ffr:gm} that
the center of the affine vertex algebra at the critical level
is closely related to Hamiltonians of the Gaudin model
describing quantum spin chain. The Hamiltonians are obtained
by the application of the Sugawara operators to the vacuum
vector of the vertex algebra $V_{\ka}(\gl_{m|n})$ with $\ka=n-m$.
Such an application yields elements of the center $\z(\wh\gl_{m|n})$
of $V_{n-m}(\gl_{m|n})$.
We thus obtain explicit formulas
for several families of elements of a commutative subalgebra of
$\U(t^{-1}\gl_{m|n}[t^{-1}])$.

Using a duality between the superalgebras
$\U(t^{-1}\gl_{m|n}[t^{-1}])$ and $\U(\gl_{m|n}[t])$ we also
obtain the corresponding families of commuting elements
in the superalgebra $\U(\gl_{m|n}[t])$. More general families
of commutative
subalgebras of $\U(\gl_{m|n}[t])$ can be brought in by applying some
automorphisms of this superalgebra parameterized by
diagonal or numerical matrices.

The connection with the Gaudin model is obtained by
considering the iterated comultiplication map
\beql{comult}
\U(\gl_{m|n}[t])\to \U(\gl_{m|n}[t])\ot\dots\ot\U(\gl_{m|n}[t])
\eeq
with $k$ copies of $\U(\gl_{m|n}[t])$ and then
taking the images of elements of a commutative subalgebra
of $\U(\gl_{m|n}[t])$ in the
tensor product of evaluation
representations of $\U(\gl_{m|n}[t])$. This yields
the higher Hamiltonians of the Gaudin model;
see Sec.~\ref{subsec:cssg} below for more details.
Such a scheme was used in \cite{t:qg}
to give simple and explicit
determinant-type formulas for the higher Gaudin
Hamiltonians in the case of $\gl_n$;
see also \cite{cm:ho},
\cite{ct:qs}, \cite{mtv:be}, \cite{mtv:ba}
and references therein.
The results of \cite{cm:ho} include a calculation
of the eigenvalues of the Sugawara operators
in the Wakimoto modules which we believe can be extended
to the super case as well with the use of the results of
\cite{ik:wm}.

\medskip
\noindent
{\bf 1.5. Singular vectors in Verma modules.}\quad
Since the Sugawara operators belong to the center
of the local completion of the universal enveloping
algebra $\U(\wh\gl_{m|n})$ at the critical level,
they commute, in particular,
with the elements of the Lie superalgebra $\wh\gl_{m|n}$.
Therefore, the Sugawara operators form a commuting
family of $\wh\gl_{m|n}$-endomorphisms of any Verma module
$M(\la)$ at the critical level. This means that
singular vectors of the Verma module can be constructed
by the applications of the Sugawara operators
to the highest vector. This approach goes back to
the works of Goodman and Wallach~\cite{gw:ho}
and Hayashi~\cite{h:so} where it was applied to the
classical Lie algebras and used to prove the Kac--Kazhdan
conjecture on the character of the irreducible
quotient of $M(\la)$.

One of the main results of a recent paper of
Gorelik~\cite[Theorem~1.1]{g:gv} describes the structure
of the algebra of singular vectors of generic Verma modules
at the critical level over symmetrizable affine Lie
superalgebras. Relying on this general theorem we construct
two families of algebraically
independent generators of the algebra of singular vectors
of the generic Verma module $M(\la)$ at the critical level
over the affine Lie superalgebra
$\gl_{m|n}[t,t^{-1}]\oplus\CC K\oplus \CC d$
(Theorem~\ref{thm:gensv} and Corollary~\ref{cor:gencom}).

\medskip
Since the first version of this paper appeared in the arXiv in 2009,
some related results concerning the center of the affine vertex algebra
at the critical level (the Feigin--Frenkel center)
associated with the simple Lie algebras of types $B$, $C$ and $D$
have been obtained. In particular, explicit generators of the center were constructed
in \cite{m:ff} and their Harish-Chandra images (the eigenvalues in Wakimoto modules)
were calculated in \cite{mm:yc}.

\medskip
We are grateful to Lucy Gow for help with calculations
in the case $m=n$. E.~R. was at
the University of Sydney during the completion of this work,
he wishes to warmly
thank the School of Mathematics and Statistics for hospitality.
We acknowledge the financial
support of the Australian Research Council.

\section{Manin matrices over superalgebras}\label{sec:ms}
\setcounter{equation}{0}

\subsection{MacMahon Master Theorem}
\label{subsec:mmmt}

As defined in the Introduction, the right quantum
superalgebra $\Mc_{m|n}$ is generated by the elements
$z_{ij}$ with $1\leqslant i,j\leqslant m+n$
such that the parity of $z_{ij}$ is $\bi+\bj$.
The defining relations are given in \eqref{srq}.
In what follows we will use the matrix form \eqref{mdefrel}
of these relations.
In order to explain this notation
in more detail, consider the superalgebra
\beql{tenprk}
\End\CC^{m|n}\ot\dots\ot\End\CC^{m|n}\ot\Mc_{m|n}
\eeq
with $k$ copies of $\End\CC^{m|n}$. For each $a\in\{1,\dots,k\}$
the element $Z_a$ of the superalgebra \eqref{tenprk}
is defined by the formula
\ben
Z_a=\sum_{i,j=1}^{m+n}
1^{\ot(a-1)}\ot e_{ij}\ot 1^{\ot(k-a)}\ot z_{ij}
\ts(-1)^{\bi\tss\bj+\bj}.
\een

Using the natural action of $\Sym_k$ on $(\CC^{m|n})^{\ot k}$
we represent any permutation $\si\in \Sym_k$
as an element $P_{\si}$ of the superalgebra \eqref{tenprk} with
the identity component in $\Mc_{m|n}$.
In particular, the transposition $(a\ts b)$ with $a<b$
corresponds to the element
\beql{ars}
P_{ab}=\sum_{i,j=1}^{m+n}
1^{\ot(a-1)}\ot e_{ij}\ot 1^{\ot(b-a-1)}
\ot e_{ji}\ot 1^{\ot(k-b)}\ot 1\ts (-1)^{\bj},
\eeq
which allows one to determine $P_{\si}$
by writing an arbitrary $\si\in\Sym_k$ as a product
of transpositions.
Recall also that if
$x,x'$ are homogeneous elements of a superalgebra $\Ac$ and $y,y'$
are homogeneous elements of a superalgebra $\Bc$ then the
product in the superalgebra $\Ac\ot\Bc$ is calculated
with the use of the sign rule
\ben
(x\ot y)(x'\ot y')=(xx'\ot yy')\ts(-1)^{\deg{y}\ts\deg{x'}}.
\een
Note that all the elements $Z_a$ and $P_{ab}$ are even
and $P_{ab}$ commutes with $Z_c$ if $c\ne a,b$.
Moreover, $P_{ab}Z_a=Z_bP_{ab}$.

For each $a=1,\dots,k$ the supertrace $\str_a$
with respect to the $a$-th copy of $\End\CC^{m|n}$
in \eqref{tenprk} is the linear map
\ben
\str_a: \big(\End\CC^{m|n}\big)^{\ot\ts k}\ot\Mc_{m|n}
\to \big(\End\CC^{m|n}\big)^{\ot\ts k-1}\ot\Mc_{m|n},
\een
defined by
\ben
\bal
\str_a:x_1\ot\dots \ot x_{a-1}&\ot e_{ij}\ot
x_{a+1}\ot \dots\ot x_k\ot y\\
{}&{}\mapsto \de_{ij}\ts x_1\ot\dots \ot x_{a-1}\ot x_{a+1}
\ot \dots\ot x_k\ot y\ts(-1)^{\bi}.
\eal
\een
In the case $k=1$ this definition clearly agrees with
\eqref{supertr}. The subscripts like $\str^{}_{1,\dots,k}$ of
the supertrace will
indicate that it is taken over the copies $1,2,\dots,k$
of the superalgebra $\End\CC^{m|n}$
in \eqref{tenprk}.
The following cyclic property of the supertrace will often
be used: if $X=[x_{ij}]$ and $Y=[y_{ij}]$ are even matrices
with pairwise super-commuting entries, $[x_{ij},y_{kl}]=0$,
then
\ben
\str(XY-YX)=0.
\een

We will need the properties of Manin matrices
given in the next proposition. In the even case ($n=0$)
they are formulated in \cite[Proposition~18]{cfr:ap}.
Recall that $H_k$ and $A_k$ denote the respective images
of the symmetrizer and antisymmetrizer
in \eqref{tenprk}; see \eqref{symm}.

\bpr\label{prop:ansym}
We have the identities in the superalgebra \eqref{tenprk},
\beql{ak}
A_k Z_1\dots Z_k A_k=A_k Z_1\dots Z_k
\eeq
and
\beql{hk}
H_k Z_1\dots Z_k H_k= Z_1\dots Z_k H_k.
\eeq
\epr

\bpf
To prove \eqref{ak} it suffices to show that for any $a=1,\dots,k-1$
we have
\beql{ampam}
A_k Z_1\dots Z_k\tss P_{a\ts a+1}=-A_k Z_1\dots Z_k.
\eeq
Since $A_k=A_k\ts (1-P_{a\ts a+1})/2$,
it is enough to consider the case $k=2$.
The relation \eqref{ampam} then reads
\ben
\frac{1-P_{12}}{2}\ts Z_1\tss Z_2\tss P_{12}=-\frac{1-P_{12}}{2}\ts Z_1\tss Z_2
\een
which an equivalent form of \eqref{mdefrel}.
Similarly, the proof of \eqref{hk} reduces to checking that
\ben
P_{a\ts a+1}\tss Z_1\dots Z_k\tss H_k= Z_1\dots Z_k H_k.
\een
This follows again from \eqref{mdefrel} written in the form
\ben
P_{12}\ts Z_1\tss Z_2\tss \frac{1+P_{12}}{2}=Z_1\tss Z_2\tss \frac{1+P_{12}}{2}.
\vspace{-1.5em}
\een
\epf

We are now in a position to prove the MacMahon Master Theorem
for the right quantum superalgebra $\Mc_{m|n}$.
We use the notation \eqref{bf}.

\bth\label{thm:mmmt}
We have the identity
\ben
{\rm Bos}\times {\rm Ferm}=1.
\een
\eth

\bpf
It is sufficient to show that for any integer $k\geqslant 1$
we have the identity in the superalgebra \eqref{tenprk}
\ben
\sum_{r=0}^k (-1)^{k-r}\str_{1,\dots,\ts r}H_rZ_1\dots Z_r\times
\str_{r+1,\dots,\ts k}A_{\{r+1,\dots,\ts k\}}Z_{r+1}\dots Z_k=0,
\een
where $A_{\{r+1,\dots,\ts k\}}$ denotes the antisymmetrizer
in \eqref{tenprk}
over the copies of $\End\CC^{m|n}$ labeled by $r+1,\dots,k$.
The identity can be written as
\beql{idenah}
\sum_{r=0}^k (-1)^r\str_{1,\dots,\ts k}H_rA_{\{r+1,\dots,\ts k\}}
Z_1\dots Z_k=0.
\eeq
Our next step is to show that the product
of the symmetrizer and antisymmetrizer in \eqref{idenah}
can be replaced as follows:
\beql{repla}
H_rA_{\{r+1,\dots,\ts k\}}\mapsto\frac{r(k-r+1)}{k}\ts
H_rA_{\{r,\dots,\ts k\}}+\frac{(r+1)(k-r)}{k}\ts
H_{r+1}A_{\{r+1,\dots,\ts k\}}.
\eeq
Indeed, the right hand side of \eqref{repla} equals
\ben
\bal
\frac{r(k-r+1)}{k}\ts
H_r&\Big(\frac{1}{k-r+1}\ts A_{\{r+1,\dots,\ts k\}}-
\frac{k-r}{k-r+1}\ts
A_{\{r+1,\dots,\ts k\}}P_{r,r+1}A_{\{r+1,\dots,\ts k\}}\Big)\\
{}+\frac{(r+1)(k-r)}{k}\ts
&\Big(\frac{1}{r+1}\ts H_r+
\frac{r}{r+1}\ts
H_rP_{r,r+1}H_r\Big)
A_{\{r+1,\dots,\ts k\}}.
\eal
\een
Since $H_r$ commutes with $A_{\{r+1,\dots,\ts k\}}$,
using the cyclic property of the supertrace we get
\begin{multline}
\str_{1,\dots,\ts k}H_rA_{\{r+1,\dots,\ts k\}}
P_{r,r+1}A_{\{r+1,\dots,\ts k\}}
Z_1\dots Z_k\\
{}=\str_{1,\dots,\ts k}H_r
P_{r,r+1}A_{\{r+1,\dots,\ts k\}}
Z_1\dots Z_kA_{\{r+1,\dots,\ts k\}}.
\non
\end{multline}
Now apply the first relation of Proposition~\ref{prop:ansym}
to write
this element as
\ben
\str_{1,\dots,\ts k}H_r
P_{r,r+1}A_{\{r+1,\dots,\ts k\}}
Z_1\dots Z_k.
\een
Similarly, using the second relation of
Proposition~\ref{prop:ansym}
and the cyclic property of the supertrace, we get
\ben
\str_{1,\dots,\ts k}H_rP_{r,r+1}H_r
A_{\{r+1,\dots,\ts k\}}
Z_1\dots Z_k=\str_{1,\dots,\ts k}H_r
P_{r,r+1}A_{\{r+1,\dots,\ts k\}}
Z_1\dots Z_k
\een
thus showing that the left hand side of \eqref{idenah}
remains unchanged after the replacement \eqref{repla}.
Then the telescoping sum in \eqref{idenah}
vanishes and the proof is complete.
\epf

We now give explicit formulas for the summands in
\eqref{bf} in terms of the generators $z_{ij}$
of the superalgebra $\Mc_{m|n}$.
We fix some notation first. Given any permutation $\si\in\Sym_k$
consider the product of the elements $P_{\si}$ and
$e_{i_1j_1}\ot\dots\ot e_{i_kj_k}\ot 1$
of the superalgebra \eqref{tenprk}.
We define the coefficient $\vp(\si,I,J)$ by
writing this product
in the form
\beql{phinot}
P_{\si}\big(e_{i_1j_1}\ot\dots\ot e_{i_kj_k}\ot 1\big)
=\vp(\si,I,J)\cdot
e_{i_{\si^{-1}(1)}j_1}\ot\dots\ot e_{i_{\si^{-1}(k)}j_k}\ot 1,
\eeq
where $I=\{i_1,\dots,i_k\}$ and $J=\{j_1,\dots,j_k\}$
so that $\vp(\si,I,J)$ equals $1$ or $-1$
depending on the multisets $I,J$ and the permutation $\si$.
Similarly, define the sign $\psi(\si,I,J)$ by the relation
\ben
\big(e_{i_1j_1}\ot\dots\ot e_{i_kj_k}\ot 1\big)\tss P_{\si^{-1}}
=\psi(\si,I,J)\cdot
e_{i_1j_{\si^{-1}(1)}}\ot\dots\ot e_{i_kj_{\si^{-1}(k)}}\ot 1.
\een
We also set
\beql{ganot}
\ga(I,J)=\sum_{a}\bi_a\tss\bj_a+
\sum_{a<b}(\bi_a+\bj_a)(\bi_b+\bj_b)
\eeq
and let $\al_i$ denote
the multiplicity of index $i$ in the given multiset $I$.

\bpr\label{prop:extsh}
We have the relation
\ben
\bal
\str_{1,\dots,\ts k}\ts A_k Z_1\dots Z_k
{}&=\sum_{I}\frac{1}{\al_{m+1}!\dots\al_{m+n}!}\\
{}&\times{}\sum_{\si\in\Sym_k}\sgn\si\cdot\vp(\si,\si I,I)\cdot
z_{i_{\si(1)}i_1}\dots z_{i_{\si(k)}i_k}(-1)^{\ga(\si I,I)},
\eal
\een
summed over multisets
$I=\{i_1\geqslant\dots\geqslant i_l\geqslant m+1>i_{l+1}>\dots>i_k\}$
with $l=0,\dots,k$,
where $\si I=\{i_{\si(1)},\dots,i_{\si(k)}\}$.
Moreover,
\ben
\bal
\str_{1,\dots,\ts k}\ts  Z_1\dots Z_kH_k
{}&=\sum_{I}\frac{1}{\al_{1}!\dots\al_{m}!}\\
{}&\times{}\sum_{\si\in\Sym_k}\psi(\si, I,\si I)\cdot
z_{i_1i_{\si(1)}}\dots z_{i_ki_{\si(k)}}
(-1)^{\ga(I,\si I)},
\eal
\een
summed over multisets
$I=\{i_1\leqslant\dots\leqslant i_l\leqslant m<i_{l+1}<\dots<i_k\}$
with $l=0,\dots,k$.
\epr

\bpf
Write
\beql{expa}
A_k Z_1\dots Z_k=\sum_{I,J}
e_{i_1j_1}\ot\dots\ot e_{i_kj_k}\ot
Z^{\ts i_1\dots\tss i_k}_{\ts j_1\dots\tss j_k},
\eeq
summed over multisets
$I=\{i_1,\dots,i_k\}$ and $J=\{j_1,\dots,j_k\}$,
where $Z^{\ts i_1\dots\tss i_k}_{\ts j_1\dots\tss j_k}\in\Mc_{m|n}$.
By Proposition~\ref{prop:ansym},
for each $a=1,\dots,k-1$ we have
\ben
P_{a,a+1}A_k Z_1\dots Z_k=A_k Z_1\dots Z_kP_{a,a+1}
=-A_k Z_1\dots Z_k.
\een
Therefore, the coefficients in the expansion \eqref{expa}
satisfy
\ben
Z^{\ts i_1\dots\tss i_{a+1}i_a\dots\tss i_k}_{\ts j_1\dots\tss j_k}
=-Z^{\ts i_1\dots\tss i_k}_{\ts j_1\dots\tss j_k}
(-1)^{\bi_{a}\bi_{a+1}+\bi_{a}\bj_{a}+\bi_{a+1}\bj_{a}}
\een
and
\ben
Z^{\ts i_1\dots\tss i_k}_{\ts j_1\dots\tss j_{a+1}j_a\dots\tss j_k}
=-Z^{\ts i_1\dots\tss i_k}_{\ts j_1\dots\tss j_k}
(-1)^{\bj_{a}\bj_{a+1}+\bj_{a}\bi_{a+1}+\bj_{a+1}\bi_{a+1}}.
\een
Either of these relations implies that if $i_a=i_{a+1}$ then
\beql{alte}
Z^{\ts i_1\dots\tss i_{a}i_{a+1}\dots\tss i_k}_{\ts i_1
\dots\tss i_{a}i_{a+1}\dots\tss i_k}
=-Z^{\ts i_1\dots\tss i_{a}i_{a+1}\dots\tss i_k}_{\ts i_1
\dots\tss i_{a}i_{a+1}\dots\tss i_k}
(-1)^{\bi_{a}}
\eeq
so that the coefficient vanishes if $\bi_{a}=0$.
Moreover, for any multiset $I$ we have
\beql{equa}
Z^{\ts i_1\dots\tss i_{a+1}i_a\dots\tss i_k}_{\ts i_1
\dots\tss i_{a+1}i_a\dots\tss i_k}=
Z^{\ts i_1\dots\tss i_k}_{\ts i_1\dots\tss i_k}.
\eeq
By the definition of the supertrace, we have
\ben
\str_{1,\dots,\ts k}\ts A_k Z_1\dots Z_k
=\sum_I Z^{\ts i_1\dots\tss i_k}_{\ts i_1\dots\tss i_k}
(-1)^{\bi_1+\dots+\bi_k}.
\een
Using \eqref{alte} and \eqref{equa}, we can rewrite
this expression by taking the summation over multisets
of the form
$I=\{i_1\geqslant\dots\geqslant i_l\geqslant m+1>i_{l+1}>\dots>i_k\}$
with $l=0,\dots,k$ so that
\ben
\str_{1,\dots,\ts k}\ts A_k Z_1\dots Z_k
=\sum_I \frac{k!}{\al_{m+1}!\dots\al_{m+n}!}\ts
Z^{\ts i_1\dots\tss i_k}_{\ts i_1\dots\tss i_k}
(-1)^{\bi_1+\dots+\bi_k}.
\een

On the other hand, multiplying the elements
of the superalgebra \eqref{tenprk} we find that
\ben
\bal
A_k Z_1\dots Z_k=\frac{1}{k!}&\ts\sum_{\si\in\Sym_k}
\sgn\si\cdot P_{\si}\\
{}\times{}&\sum_{I,J}
e_{i_1j_1}\ot\dots\ot e_{i_kj_k}\ot z_{i_1j_1}\dots z_{i_kj_k}
(-1)^{\sum_a (\bi_a\tss\bj_a+\bj_a)
+\sum_{a<b}(\bi_a+\bj_a)(\bi_b+\bj_b)}.
\eal
\een
Then using the notation \eqref{phinot} and \eqref{ganot}
we can write
\ben
Z^{\ts i_1\dots\tss i_k}_{\ts j_1\dots\tss j_k}
=\frac{1}{k!}\ts\sum_{\si\in\Sym_k}\sgn\si\cdot\vp(\si,\si I,J)\cdot
z_{i_{\si(1)}j_1}\dots z_{i_{\si(k)}j_k}
(-1)^{\ga(\si I,J)}\ts (-1)^{\bj_1+\dots+\bj_k},
\een
which completes the proof of the first relation.
The second relation is verified by a similar
argument with the use of the identities
\ben
P_{a,a+1}Z_1\dots Z_kH_k =Z_1\dots Z_kH_kP_{a,a+1}
=Z_1\dots Z_kH_k
\een
implied by Proposition~\ref{prop:ansym}.
\epf

\bre\label{rem:orde}
(i)\quad Relation \eqref{equa} can be used to get alternative
formulas for the elements $\str_{1,\dots,\ts k}\ts A_k Z_1\dots Z_k$
and $\str_{1,\dots,\ts k}\ts  Z_1\dots Z_kH_k$ with
the summation over multisets of the form
$\{i_1<\dots< i_l< m+1\leqslant i_{l+1}\leqslant\dots\leqslant i_k\}$
and $\{i_1>\dots> i_l> m
\geqslant i_{l+1}\geqslant\dots\geqslant i_k\}$,
respectively.
\par
(ii)\quad
Taking the particular cases $n=0$ and $m=0$
in Theorem~\ref{thm:mmmt} we thus get
a new proof of the MacMahon Master Theorem
for the right quantum algebra; cf.
\cite{fh:np}, \cite{glz:qm}, \cite{kp:ne}.
Some other noncommutative versions
of the theorem associated with
quantum algebras or classical Lie algebras
can be found e.g. in \cite{gps:sp} and \cite[Ch.~7]{m:yc};
see also references therein.
\qed
\ere

\subsection{Affine right quantum superalgebras}
\label{subsec:arqs}

In accordance to the definition we gave in the Introduction,
the affine right quantum superalgebra $\wh\Mc_{m|n}$
is generated by the elements $z^{(r)}_{ij}$
of parity $\bi+\bj$, where $r$ runs over the set of positive
integers and $1\leqslant i,j\leqslant m+n$.
The defining relations \eqref{mdefrelaff} can be written
more explicitly as follows: for all positive integers $p$
we have
\beql{srqrel}
\sum_{r+s=p}\Big([z^{(r)}_{ij}, z^{(s)}_{kl}]
-[z^{(r)}_{kj}, z^{(s)}_{il}](-1)^{\bi\bj+\bi\bk+\bj\bk}\Big)=0,
\eeq
summed over nonnegative integers $r$ and $s$,
where we set $z^{(0)}_{ij}=\de_{ij}$. Equivalently,
in terms of the formal power series \eqref{serz} they take the form
\beql{drseries}
[z_{ij}(u), z_{kl}(u)]=[z_{kj}(u), z_{il}(u)]
(-1)^{\bi\bj+\bi\bk+\bj\bk},
\qquad i,j,k,l\in \{1,\dots,m+n\}.
\eeq
We will keep
the notation $z^{\tss\prime}_{ij}(u)$ for the entries
of the inverse matrix $Z^{-1}(u)$.

\bpr\label{prop:ev}
The mapping
\beql{eval}
z_{ij}(u)\mapsto \delta_{ij}+z_{ij}\tss u
\eeq
defines a surjective homomorphism $\wh\Mc_{m|n}\to\Mc_{m|n}$. Moreover,
the mapping
\beql{embed}
z_{ij}\mapsto z_{ij}^{(1)}
\eeq
defines an embedding $\Mc_{m|n}\hookrightarrow \wh\Mc_{m|n}$.
\epr

\bpf
This follows easily from the defining relations of
$\Mc_{m|n}$ and $\wh\Mc_{m|n}$.
The injectivity of the map \eqref{embed}
follows from the observation that the composition of \eqref{embed}
and \eqref{eval} yields the identity map on $\Mc_{m|n}$.
\epf

We will identify $\Mc_{m|n}$ with a subalgebra of $\wh\Mc_{m|n}$
via the embedding \eqref{embed}.

\bde\label{def:ber}
The {\it Berezinian\/} of the matrix $Z(u)$ is a formal
power series in $u$ with coefficients in $\wh\Mc_{m|n}$
defined by the formula
\begin{align}
\Ber Z(u){}&=\sum_{\si\in\Sym_m}\sgn\si\cdot
z_{\si(1)1}(u)\dots z_{\si(m)m}(u)
\non\\
\label{berdef}
{}&\times{}
\sum_{\tau\in\Sym_n}\sgn\tau\cdot
z^{\tss\prime}_{m+1,m+\tau(1)}(u)\dots
z^{\tss\prime}_{m+n,m+\tau(n)}(u).
\end{align}
The image of the Berezinian $\Ber Z(u)$ under the homomorphism
\eqref{eval} will be denoted by $\Ber(1+uZ)$. This is a formal
power series in $u$ with coefficients in the right quantum
superalgebra $\Mc_{m|n}$.
\qed
\ede

An alternative formula for the Berezinian $\Ber Z(u)$
will be given in Corollary~\ref{cor:altdef} below.

Our goal now is to derive a quasideterminant factorization
of $\Ber Z(u)$ and then use it in the proof of the identities
\eqref{charferm}--\eqref{newton}. We start by providing
some symmetries of the superalgebra $\wh\Mc_{m|n}$.

\ble\label{lem:manpr}
If $Z$ is a Manin matrix, then for any
positive integer $r$ the following identity holds:
\beql{manpr}
(1-P_{12})\sum_{k+l=r}[Z_1^{\tss k},Z_2^{\tss l}]=0,
\eeq
where $k$ and $l$ run over nonnegative integers.
\ele

\bpf
We use induction on $r$.
The identity is trivial for $r=1$, while for $r=2$
it is equivalent to the definition of a Manin matrix;
see \eqref{mdefrel}. Suppose that $r\geqslant 3$.
Note first that
using the induction hypothesis we get
\ben
\bal
&(1-P_{12})Z_2Z_1\sum_{k+l=r-2}[Z_1^k,Z_2^{\tss l}]
=(1-P_{12})Z_1Z_2\sum_{k+l=r-2}[Z_1^k,Z_2^{\tss l}]\\
{}={}&{}(1-P_{12})Z_1Z_2P_{12}\sum_{k+l=r-2}[Z_1^k,Z_2^{\tss l}]
=-(1-P_{12})Z_2Z_1\sum_{k+l=r-2}[Z_1^k,Z_2^{\tss l}]
\eal
\een
and so
\beql{vansum}
(1-P_{12})Z_2Z_1\sum_{k+l=r}[Z_1^{\tss k-1},Z_2^{\tss l-1}]=0.
\eeq
Now for $k\geqslant 1$ write
\ben
[Z_1^{\tss k},Z_2^{\tss l}]=[Z_1,Z^{\tss l}_2]\ts
Z_1^{\tss k-1}
+Z_1[Z_1^{\tss k-1},Z_2^{\tss l}].
\een
If $l\geqslant 1$ then we also have
\ben
[Z_1,Z^{\tss l}_2]\ts
Z_1^{\tss k-1}=[Z_1,Z_2]\ts Z_2^{\tss l-1}
Z_1^{\tss k-1}+Z_2[Z_1,Z^{\tss l-1}_2]\ts Z_1^{\tss k-1}.
\een
Therefore, due to \eqref{mdefrel} we obtain
\ben
(1-P_{12})\sum_{k+l=r}[Z_1^{\tss k},Z_2^{\tss l}]
=(1-P_{12})\sum_{k+l=r}
\Big(Z_2[Z_1,Z^{\tss l-1}_2]\ts Z_1^{\tss k-1}
+Z_1[Z_1^{\tss k-1},Z_2^{\tss l}]\Big),
\een
where $k$ and $l$ run over positive integers.
By the induction hypothesis,
\ben
\bal
(1-P_{12})\sum_{k+l=r}Z_1[Z_1^{\tss k-1},Z_2^{\tss l}]
{}&=(1-P_{12})Z_1P_{12}\sum_{k+l=r}[Z_1^{\tss k-1},Z_2^{\tss l}]\\
{}&=-(1-P_{12})\sum_{k+l=r}Z_2[Z_1^{\tss k-1},Z_2^{\tss l}].
\eal
\een
Now for $k\geqslant 2$ write
\ben
[Z_1^{\tss k-1},Z_2^{\tss l}]=[Z_1,Z_2^{\tss l}]\tss Z_1^{\tss k-2}
+Z_1[Z_1^{\tss k-2},Z_2^{\tss l}].
\een
Hence, taking appropriate restrictions on
the summation indices $k$ and $l$ we get
\ben
\bal
&(1-P_{12})\sum_{k+l=r}[Z_1^{\tss k},Z_2^{\tss l}]\\
{}={}&(1-P_{12})\sum_{k+l=r}
\Big(Z_2[Z_1,Z^{\tss l-1}_2]\ts Z_1^{\tss k-1}
-Z_2[Z_1,Z_2^{\tss l}]\tss Z_1^{\tss k-2}
-Z_2Z_1[Z_1^{\tss k-2},Z_2^{\tss l}]\Big)
\eal
\een
which is zero by \eqref{vansum}.
\epf

\bpr\label{prop:inv}
The mapping
\beql{omegan}
\om:Z(u)\mapsto Z^{-1}(-u)
\eeq
defines an involutive automorphism of the
superalgebra $\wh\Mc_{m|n}$ which is identical
on the subalgebra $\Mc_{m|n}$.
\epr

\bpf
Set $\check Z=1-Z(-u)$ so that $\check Z$ is a formal power series
in $u$ without a constant term. Then $Z^{-1}(-u)$ can be written
in the form
\ben
Z^{-1}(-u)=\sum_{k=0}^{\infty} \check Z^{\tss k}.
\een
Note that the right hand side is a well-defined power series
in $u$ with coefficients in $\wh\Mc_{m|n}$.
By Lemma~\ref{lem:manpr}, for any positive integer $r$ we have
\ben
(1-P_{12})\sum_{k+l=r}[\check Z_1^{\tss k},\check Z_2^{\tss l}]=0,
\een
which implies
\ben
(1-P_{12})[Z_1^{-1}(-u),Z_2^{-1}(-u)]=0
\een
thus proving that $\om$ is a homomorphism.

Applying $\om$ to both sides of
$
\om\big(Z(u)\big)\ts Z(-u)=1
$
we get $\om^2\big(Z(u)\big)\ts Z^{-1}(u)=1$,
so that $\om$ is an involutive automorphism of $\wh\Mc_{m|n}$.
The second statement is clear.
\epf

We also introduce the {\it affine
left quantum superalgebra\/} $\wh\Mc^{\circ}_{m|n}$.
It is generated by a countable set of elements $y^{(r)}_{ij}$
of parity $\bi+\bj$, where $r$ runs over the set of positive
integers. The defining relations of $\wh\Mc^{\circ}_{m|n}$
take the form
\beql{mdefrelaffl}
[Y_1(u),Y_2(u)]\ts(1-P_{12})=0,
\eeq
where $u$ is a formal variable, and the matrix elements
of the matrix $Y(u)=[y_{ij}(u)]$ are the formal power series
\beql{serzl}
y_{ij}(u)=\de_{ij}+y^{(1)}_{ij}\ts u+y^{(2)}_{ij}\ts u^2+\cdots.
\eeq
The defining relations \eqref{mdefrelaffl} can be written
as
\beql{drseriesl}
[y_{ij}(u), y_{kl}(u)]+[y_{kj}(u), y_{il}(u)]
(-1)^{\bi\bj+\bi\bk+\bj\bk}=0,
\qquad i,j,k,l\in \{1,\dots,m+n\},
\eeq
or, equivalently,
for all positive integers $p$
we have
\beql{srqrell}
\sum_{r+s=p}\Big([y^{(r)}_{ij}, y^{(s)}_{kl}]
+[y^{(r)}_{kj}, y^{(s)}_{il}](-1)^{\bi\bj+\bi\bk+\bj\bk}\Big)=0,
\eeq
summed over nonnegative integers $r$ and $s$,
where we set $y^{(0)}_{ij}=\de_{ij}$.

It is straightforward to verify that the superalgebra
$\wh\Mc^{\circ}_{m|n}$ is isomorphic to the affine right
quantum superalgebra $\wh\Mc_{m|n}$. An isomorphism
can be given by the supertransposition map
$Y(u)\mapsto Z(u)^t$ so that
$y_{ij}(u)\mapsto z_{ji}(u)(-1)^{\bi\tss\bj+\bi}$.

\bpr\label{prop:isomlr}
The mapping
\beql{isomze}
\ze:y_{ij}(u)\mapsto z^{\tss\prime}_{m+n-i+1,m+n-j+1}(u)
\eeq
defines an isomorphism $\wh\Mc^{\circ}_{n|m}\to \wh\Mc_{m|n}$.
\epr

\bpf
Observe that the mapping
$y_{ij}(u)\mapsto z_{m+n-i+1,m+n-j+1}(-u)$
defines an isomorphism $\wh\Mc^{\circ}_{n|m}\to \wh\Mc_{m|n}$.
This follows easily from the defining relations of
the affine left and right quantum superalgebras.
It remains to note that $\ze$ is the composition
of this isomorphism and the automorphism $\om$ defined
in Proposition~\ref{prop:inv}.
\epf

We will now adapt the arguments used by Gow~\cite{g:rl, g:gd}
for the Yangian of the Lie superalgebra $\gl_{m|n}$
to derive a quasideterminant decomposition
of the Berezinian $\Ber Z(u)$; see Definition~\ref{def:ber}.

If $A=[a_{ij}]$ is a square matrix over a ring with $1$,
then its $ij$-{\it th quasideterminant\/} is defined
if $A$ is invertible and the $ji$-th entry $(A^{-1})_{ji}$
is an invertible element of the ring.
The $ij$-th quasideterminant is then given
by
\ben
\left|\begin{matrix}a_{11}&\dots&a_{1j}&\dots&a_{1N}\\
                                   &\dots&      &\dots&      \\
                             a_{i1}&\dots&\boxed{a_{ij}}&\dots&a_{iN}\\
                                   &\dots&      &\dots&      \\
                             a_{N1}&\dots&a_{Nj}&\dots&a_{NN}
                \end{matrix}\right|=\big((A^{-1})_{ji}\big)^{-1}.
\een
We refer the reader to \cite{ggrw:q}, \cite{gkllrt:ns}
and references therein for the properties
and applications of the quasideterminants.

We will need the Gauss decompositions of the matrices
$Z(u)$ and $Y(u)$.
There exist unique matrices $D(u)$, $E(u)$ and $F(u)$
whose entries are formal power series in $u$ with coefficients
in $\wh\Mc_{m|n}$
such that
$
D(u)=\text{\rm diag\ts}\big[d_1(u),\dots,d_{m+n}(u)\big]
$
and
\ben
E(u)=\left[
\begin{matrix}
1&e_{12}(u)&\dots&e_{1,m+n}(u)\\
0&1&\dots&e_{2,m+n}(u)\\
\vdots&\vdots&\ddots&\vdots\\
0&0&\dots&1
\end{matrix}
\right],
\qquad
F(u)=\left[
\begin{matrix}
1&0&\dots&0\\
f_{21}(u)&1&\dots&0\\
\vdots&\vdots&\ddots&\vdots\\
f_{m+n,1}(u)&f_{m+n,2}(u)&\dots&1
\end{matrix}
\right]
\een
satisfying $Z(u)=F(u)\ts D(u)\ts E(u)$.
Explicit formulas for the entries of the matrices
$D(u)$, $E(u)$ and $F(u)$ can be given
in terms of quasideterminants; see \cite{ggrw:q}.
In particular,
\beql{hmqua}
d_i(u)=\left|\begin{matrix} z_{11}(u)&\dots&z_{1,i-1}(u)&z_{1i}(u)\\
                          \vdots&\ddots&\vdots&\vdots\\
        z_{i-1,1}(u)&\dots&z_{i-1,i-1}(u)&z_{i-1,i}(u)\\
        z_{i1}(u)&\dots&z_{i,i-1}(u)&\boxed{z_{ii}(u)}\\
           \end{matrix}\right|,
\eeq
for $i=1,\dots,m+n$.

We write the Gauss decomposition
of the matrix $Y(u)$ as
$Y(u)=F^{\circ}(u)\ts D^{\circ}(u)\ts E^{\circ}(u)$ so that
the corresponding entries
$d^{\tss\circ}_i(u),e^{\circ}_{ij}(u),f^{\circ}_{ij}(u)$
of the matrices $D^{\circ}(u)$,
$E^{\circ}(u)$ and $F^{\circ}(u)$
are formal power series in $u$ with coefficients
in $\wh\Mc^{\circ}_{m|n}$. These entries are found
by the same formulas as above with the $z_{ij}(u)$
respectively replaced by $y_{ij}(u)$.

\ble\label{lem:imden}
Under the isomorphism $\ze:\wh\Mc^{\circ}_{n|m}\to\wh\Mc_{m|n}$
defined in \eqref{isomze}
we have
\ben
\ze:d^{\tss\circ}_k(u)\mapsto d_{m+n-k+1}(u)^{\tss-1},\qquad
k=1,\dots,n+m.
\een
\ele

\bpf
The entries of the inverse matrix $Z^{-1}(u)$ are found
from the decomposition
$Z^{-1}(u)=E^{-1}(u)\ts D^{-1}(u)\ts F^{-1}(u)$ so that
\beql{zprex}
z^{\tss\prime}_{i^{\tss\prime}j^{\tss\prime}}(u)
=\sum_{k\tss\leqslant\ts i,j}
e'_{i^{\tss\prime}k^{\tss\prime}}(u)\tss
d_{k^{\tss\prime}}(u)^{-1}\tss
f^{\tss\prime}_{k^{\tss\prime}j^{\tss\prime}}(u),
\eeq
where $e'_{ij}(u)$ and $f'_{ij}(u)$ denote the entries
of the matrices $E^{-1}(u)$ and $F^{-1}(u)$, respectively,
and we set $i^{\tss\prime}=m+n-i+1$.
On the other hand, the entries of $Y(u)$ are found by
\beql{yprex}
y_{ij}(u)=
\sum_{k\tss\leqslant\ts i,j}f^{\circ}_{i\tss k}(u)
\tss d^{\tss\circ}_{k}(u)
\tss e^{\circ}_{kj}(u).
\eeq
Hence, by \eqref{isomze} we have
\ben
\ze(d^{\tss\circ}_1(u))=\ze(y_{11}(u))
=z^{\tss\prime}_{1'1'}(u)=d_{1'}(u)^{-1}.
\een
Now, comparing \eqref{zprex} and \eqref{yprex},
and arguing by induction we find that
\ben
\ze:f^{\circ}_{i\tss k}(u)\mapsto e'_{i^{\tss\prime}k^{\tss\prime}}(u),
\qquad
e^{\circ}_{kj}(u)\mapsto
f^{\tss\prime}_{k^{\tss\prime}j^{\tss\prime}}(u),
\qquad
d^{\tss\circ}_{k}(u)\mapsto d_{k^{\tss\prime}}(u)^{-1},
\een
as required.
\epf

\bth\label{thm:gauss}
The Berezinian $\Ber Z(u)$
admits the quasideterminant factorization
in the superalgebra $\wh\Mc_{m|n}[[u]]${\rm :}
\ben
\Ber Z(u)=d^{}_{1}(u)\dots d^{}_m(u)\ts d^{-1}_{m+1}(u)\dots
d^{-1}_{m+n}(u).
\een
\eth

\bpf
The upper-left $m\times m$ submatrix of $Z(u)$
is a matrix with even entries satisfying
\eqref{rq}. A quasideterminant decomposition
of the corresponding determinant was obtained in
\cite[Lemma~8]{cfr:ap} and it takes the form
\beql{qddec}
\sum_{\si\in\Sym_m}\sgn\si\cdot
z_{\si(1)1}(u)\dots z_{\si(m)m}(u)=
d^{}_{1}(u)\dots d^{}_m(u).
\eeq
The second factor on the right hand side of \eqref{berdef}
is the image of the determinant
\ben
\sum_{\tau\in\Sym_n}\sgn\tau\cdot
y_{n,\tau(n)}(u)\dots y_{1,\tau(1)}(u)\in \wh\Mc^{\circ}_{n|m}[[u]]
\een
under the isomorphism \eqref{isomze}. However,
the upper-left $n\times n$ submatrix of $Y(u)$
is a matrix with even entries
whose transpose satisfies \eqref{rq}.
The corresponding quasideterminant decomposition
is proved in the same way as \eqref{qddec} (see \cite{cfr:ap}),
so that
\ben
\sum_{\tau\in\Sym_n}\sgn\tau\cdot
y_{n,\tau(n)}(u)\dots y_{1,\tau(1)}(u)
=d^{\tss\circ}_n(u)\dots d^{\tss\circ}_1(u).
\een
By Lemma~\ref{lem:imden}, the image of this determinant
under $\ze$ is $d^{-1}_{m+1}(u)\dots
d^{-1}_{m+n}(u)$.
\epf

An alternative
factorization of the Berezinian involving different
quasideterminants is provided
by Corollary~\ref{cor:altdef} below.

\bre\label{rem:yang} (i)\quad
In the super-commutative specialization the quasideterminant decompositions
of the Berezinian given in Theorem~\ref{thm:gauss} and in Corollary~\ref{cor:altdef} below
turn into the decompositions originally found by I.~Gelfand and V.~Retakh;
cf. \cite[Theorem~3.8.1]{ggrw:q}.
\par
(ii)\quad
The quantum Berezinian $\Ber T(u)$ of the
generator matrix of the Yangian for the Lie superalgebra $\gl_{m|n}$
was introduced by Nazarov~\cite{n:qb}. The quasideterminant
decomposition of $\Ber T(u)$ found by Gow~\cite[Theorem~1]{g:rl} can
be obtained as a particular case of Theorem~\ref{thm:gauss} by
taking $Z(u)=e^{-\di_u}T(u)$. The latter is a Manin matrix which
follows easily from the defining relations of the Yangian; cf.
\cite{cf:mm, cfr:ap}. Hence, all identities for the Berezinian of
the matrix $Z(u)$ obtained in this paper imply the corresponding
counterparts for $\Ber T(u)$. \qed \ere

\subsection{Berezinian identities}

Consider the Berezinian $\Ber(1+uZ)$, where
$Z$ is a Manin matrix; see Definition~\ref{def:ber}.
The expression $\Ber(1+uZ)$ is a formal power series in $u$
with coefficients in the right quantum superalgebra
$\Mc_{m|n}$. The next theorem provides
some identities for the coefficients of this series,
including a noncommutative analogue
of the Newton identities \eqref{newtonthm};
cf. \cite{cf:mm}, \cite{kv:be}, \cite{s:si}.

\bth\label{thm:ident}
We have the identities
\begin{align}\label{charfermthm}
\Ber(1+uZ)&=\sum_{k=0}^{\infty} u^k\ts
\str^{}_{1,\dots,k}\ts A_k Z_1\dots Z_k,\\
\label{charbosthm}
\big[\Ber(1-uZ)\big]^{-1}&=\sum_{k=0}^{\infty}
u^k\ts\str^{}_{1,\dots,k}\ts H_k Z_1\dots Z_k,\\
\label{newtonthm}
\big[\Ber(1+uZ)\big]^{-1}\ts\di_u\tss \Ber(1+uZ)
&=\sum_{k=0}^{\infty} (-u)^k\ts\str\ts Z^{k+1}.
\end{align}
\eth

\bpf
Due to the MacMahon Master Theorem (Theorem~\ref{thm:mmmt}),
identities \eqref{charfermthm} and \eqref{charbosthm}
are equivalent. Moreover, \eqref{charfermthm} is clear
for $n=0$ as the Berezinian turns into a determinant.
We will be proving \eqref{charbosthm} by induction on
$n$, assuming that $n\geqslant 1$. Let $\wt Z$ be the matrix
obtained from $Z$ by deleting the row and column $m+n$.
Set
\ben
h_k(Z)=\str^{}_{1,\dots,k}\ts H_k Z_1\dots Z_k
=\str^{}_{1,\dots,k}\ts Z_1\dots Z_kH_k;
\een
the second equality holds by the cyclic property
of the supertrace.
Applying Theorem~\ref{thm:gauss} we derive that
\ben
\Ber(1-uZ)=\Ber(1-u\wt Z)\big[(1-uZ)^{-1}\big]_{m+n,m+n}
\een
and so
\ben
\big[(1-uZ)^{-1}\big]_{m+n,m+n}\big[\Ber(1-uZ)\big]^{-1}
=\big[\Ber(1-u\wt Z)\big]^{-1}.
\een
Hence, \eqref{charbosthm} will follow if we show that
\ben
\big[(1-uZ)^{-1}\big]_{m+n,m+n}\cdot
\sum_{k=0}^{\infty}
u^k\ts h_k(Z)
=\sum_{k=0}^{\infty}
u^k\ts h_k(\wt Z),
\een
or, equivalently, that for any $r\geqslant 1$ we have
\beql{powe}
\sum_{k+l=r}(Z^k)_{m+n,m+n} \cdot h_l(Z)
=h_r(\wt Z).
\eeq
In order to verify \eqref{powe} we use a relation in the
superalgebra \eqref{tenprk},
\beql{potr}
Z_1^k=\str^{}_{2,\dots,k} Z_1\dots Z_kP_{k-1,k}\dots P_{23}P_{12},
\eeq
which follows easily by induction
with the use of the relation $\str^{}_{2}Z_2P_{12}=Z_1$.
The next lemma is a Manin matrix version
of the corresponding identities obtained in \cite{iop:gc}.

\ble\label{lem:relpow}
We have the identity
\beql{hid}
\sum_{k=1}^r Z_1^k\ts h_{r-k}(Z)=r\ts
\str^{}_{2,\dots,r} Z_1\dots Z_r H_r.
\eeq
\ele

\bpf
By \eqref{potr}, the left hand side can be written as
\beql{lhsid}
\sum_{k=1}^r \str^{}_{2,\dots,r} Z_1\dots Z_r
H_{\{k+1,\dots,r\}}P_{k-1,k}\dots P_{23}P_{12}.
\eeq
Write
\ben
H_{\{k+1,\dots,r\}}=(r-k+1)\ts H_{\{k,\dots,r\}}-
(r-k)\ts H_{\{k+1,\dots,r\}}P_{k,k+1}H_{\{k+1,\dots,r\}}.
\een
By the second relation of Proposition~\ref{prop:ansym}
and the cyclic property of the supertrace, we have
\ben
\bal
&\str^{}_{2,\dots,r} Z_1\dots Z_r
H_{\{k+1,\dots,r\}}P_{k,k+1}H_{\{k+1,\dots,r\}}
P_{k-1,k}\dots P_{23}P_{12}\\
{}={}&\str^{}_{2,\dots,r} H_{\{k+1,\dots,r\}}Z_1\dots Z_r
H_{\{k+1,\dots,r\}}P_{k,k+1}
P_{k-1,k}\dots P_{23}P_{12}\\
{}={}&\str^{}_{2,\dots,r} Z_1\dots Z_r
H_{\{k+1,\dots,r\}}P_{k,k+1}\dots P_{23}P_{12}.
\eal
\een
Hence, \eqref{lhsid} takes the form of a telescoping sum
which simplifies to become the right hand side of \eqref{hid}.
\epf

Taking into account Lemma~\ref{lem:relpow}, we can represent
\eqref{powe} in the equivalent form
\beql{eqfre}
\str^{}_{1,\dots,r} Z_1\dots Z_r H_r+
r\ts\big[\tss\str^{}_{2,\dots,r} Z_1\dots Z_r H_r\big]_{m+n,m+n}
=\str^{}_{1,\dots,r} \wt Z_1\dots \wt Z_r H_r.
\eeq
Write
\ben
Z_1\dots Z_r H_r=\sum e_{i_1j_1}\ot\dots\ot e_{i_rj_r}\ot
z_{\ts j_1,\dots,j_r}^{\ts i_1,\dots,i_r},
\een
summed over all indices $i_a,j_a\in\{1,\dots,m+n\}$, where
the $z_{\ts j_1,\dots,j_r}^{\ts i_1,\dots,i_r}$ are certain
elements of the right quantum superalgebra $\Mc_{m|n}$.
Then
\beql{tracef}
\str^{}_{1,\dots,r} Z_1\dots Z_r H_r=\sum
z_{\ts i_1,\dots,i_r}^{\ts i_1,\dots,i_r}(-1)^{\bi_1+\dots+\bi_r}.
\eeq
We have the following analogue of \eqref{alte} which is verified
in the same way:
\ben
z^{\ts i_1\dots\tss i_{a}i_{a+1}\dots\tss i_k}_{\ts i_1
\dots\tss i_{a}i_{a+1}\dots\tss i_k}
=z^{\ts i_1\dots\tss i_{a}i_{a+1}\dots\tss i_k}_{\ts i_1
\dots\tss i_{a}i_{a+1}\dots\tss i_k}
(-1)^{\bi_{a}},
\een
while the counterpart of \eqref{equa} has exactly the same form.
Hence, the index
$m+n$ may occur at most once amongst
the summation indices $i_1,\dots,i_r$.
The sum in \eqref{tracef} over the indices restricted to
$i_1,\dots,i_r\in\{1,\dots,m+n-1\}$ coincides with
the expression $\str^{}_{1,\dots,r} \wt Z_1\dots \wt Z_r H_r$,
while the sum over the multisets of indices
containing $m+n$ equals
${-r}\ts\big[\tss\str^{}_{2,\dots,r}
Z_1\dots Z_r H_r\big]_{m+n,m+n}$, thus proving
\eqref{eqfre}. This completes the proof of \eqref{charbosthm}.

Now we prove the Newton identity \eqref{newtonthm}
which can be written in the equivalent form
\ben
-\di_u\big[\Ber(1+uZ)\big]^{-1}
=\sum_{k=0}^{\infty} (-u)^k\ts\str\ts Z^{k+1}
\cdot \big[\Ber(1+uZ)\big]^{-1}.
\een
Equating the coefficients of the same powers of $u$ we can
also write this as
\ben
\sum_{k=1}^r \big(\str\ts Z^k\big)\ts h_{r-k}(Z)=r\ts
h_r(Z),\qquad r=1,2,\dots.
\een
However, this relation is immediate from Lemma~\ref{lem:relpow}
by taking the supertrace $\str_1$ over the first copy
of the superalgebra $\End\CC^{m|n}$.
\epf

\bre\label{rem:newton}
(i)\quad
In the case $n=0$ we thus get a new proof of the Newton
identities for Manin matrices
based on Lemma~\ref{lem:relpow};
cf. \cite{cf:mm}. This argument essentially
follows \cite{iop:gc}.
\par
(ii)\quad
Relations of Theorem~\ref{thm:ident} imply that
the coefficients of the powers of $u$
in the expansions \eqref{charfermthm},
\eqref{charbosthm} and \eqref{newtonthm} can be regarded
as respective specializations of the noncommutative
elementary, complete and power sums symmetric functions
of the first kind; see \cite{gkllrt:ns}. It would
be interesting to find direct specializations of the other types
of the noncommutative symmetric functions, in particular,
the ribbon Schur functions.
\qed
\ere

Using Theorem~\ref{thm:ident} we can obtain alternative
expressions for the Berezinians $\Ber Z(u)$ and $\Ber(1+uZ)$
and different quasideterminant factorizations;
cf. Definition~\ref{def:ber} and Theorem~\ref{thm:gauss}.

\bco\label{cor:altdef}
The following relations hold
\begin{align}
\big[\Ber Z(u)\big]^{-1}{}&=\sum_{\si\in\Sym_m}\sgn\si\cdot
z^{\tss\prime}_{\si(1)1}(u)\dots z^{\tss\prime}_{\si(m)m}(u)
\non\\
\label{berdefalt}
{}&\times{}
\sum_{\tau\in\Sym_n}\sgn\tau\cdot
z_{m+1,m+\tau(1)}(u)\dots z_{m+n,m+\tau(n)}(u)
\end{align}
and
\beql{qualt}
\big[\Ber Z(u)\big]^{-1}=
\bar d^{\ts-1}_{1}(u)\dots \bar d^{\ts-1}_m(u)\ts \bar d^{}_{m+1}(u)\dots
\bar d^{}_{m+n}(u),
\eeq
where $\bar d_{i}(u)$ is the quasideterminant
\ben
\bar d_i(u)=\left|\begin{matrix} \boxed{z_{ii}(u)}&\dots&z_{i,m+n}(u)\\
                          \vdots&\ddots&\vdots\\
        z_{m+n,i}(u)&\dots&z_{m+n,m+n}(u)\\
           \end{matrix}\right|.
\een
In particular, the corresponding relations hold for
the Berezinian $\Ber(1+uZ)$.
\eco

\bpf
The arguments are quite similar to those used
in the proof of Theorems~\ref{thm:gauss} and \ref{thm:ident}
so we only sketch the main steps.
We show first that the right hand
side of \eqref{berdefalt}
coincides with the product of quasideterminants on
the right hand side of \eqref{qualt}.
Using the same calculation as in \cite[Lemma~8]{cfr:ap},
we obtain the following quasideterminant decomposition:
\ben
\sum_{\tau\in\Sym_n}\sgn\tau\cdot
z_{m+1,m+\tau(1)}(u)\dots z_{m+n,m+\tau(n)}(u)
=\bar d^{}_{m+1}(u)\dots
\bar d^{}_{m+n}(u).
\een
Furthermore, the determinant
\ben
\sum_{\si\in\Sym_m}\sgn\si\cdot
z^{\tss\prime}_{\si(1)1}(u)\dots z^{\tss\prime}_{\si(m)m}(u)
\een
coincides with the image of a certain determinant
of a submatrix of $Y(u)$
under the isomorphism \eqref{isomze} which leads
to the desired factorization via a natural dual analogue
of Lemma~\ref{lem:imden}.

Observe that the matrix $Z(u)$ can be written as
$Z(u)=1+\wt Z(u)$ and $\wt Z(u)$ satisfies
\eqref{mdefrelaff}.
It is therefore sufficient to verify relations \eqref{berdefalt}
and \eqref{qualt} for matrices of the form $Z(u)=1+uZ$,
where $Z$ satisfies \eqref{mdefrel}. Since we have verified
that the right hand sides of \eqref{berdefalt}
and \eqref{qualt} coincide, it is enough to show that
\eqref{qualt} holds. We will use
induction on
$m$, assuming that $m\geqslant 1$. Let $\overline Z$ be the matrix
obtained from $Z$ by deleting the row and column $1$.
Then we need to verify that
\ben
\big[\Ber(1+uZ)\big]^{-1}=\big[(1+uZ)^{-1}\big]_{11}
\big[\Ber(1+u\overline Z)\big]^{-1},
\een
or, equivalently,
\beql{eqde}
\Ber(1+uZ)\cdot \big[(1+uZ)^{-1}\big]_{11}=\Ber(1+u\overline Z).
\eeq
Set
\ben
\si_k(Z)=\str^{}_{1,\dots,k}\ts A_k Z_1\dots Z_k.
\een
Due to \eqref{charfermthm}, the relation \eqref{eqde}
will follow if we show that
for any $r\geqslant 1$
\ben
\sum_{k+l=r}\si_k(Z)\cdot (Z^l)_{11}
=\si_r(\overline Z).
\een
However, this follows by
the same argument as for the
proof of \eqref{powe}.
\epf

\section{Sugawara operators for $\wh\gl_{m|n}$}\label{sec:so}
\setcounter{equation}{0}

Consider the Lie superalgebra
$\wh\gl_{m|n}\oplus\CC\tau$
with its commutation relations \eqref{commrel}
and \eqref{tauco}. We assume that $m$ and $n$ are
nonnegative integers. Recall also that
$V=V_{\ka}(\gl_{m|n})$ is the affine vertex algebra
at the level $\ka\in\CC$. This means that
$V$ is equipped with
the additional data $(Y, D, 1)$, where $1 \in V$
is the vacuum vector,
the state-field correspondence $Y$ is a map
$$
Y: V \rightarrow \End V[[z,z^{-1}]],
$$
the infinitesimal translation $D$ is an operator
$D: V \rightarrow V$. These data
satisfy the vertex algebra axioms; see e.g.
\cite{fb:va}, \cite{k:va}.
For $a \in V$ we write
\ben
Y(a,z) = \sum_{r\in\ZZ} a_{(r)} z^{-r-1}, \quad
a_{(r)} \in \End V.
\een
In particular, for all $a,b \in V$ we have
$a_{(r)}\ts b = 0$ for $r \gg 0$.
The span in $\End V$ of all {\it Fourier coefficients\/} $a_{(r)}$
of all vertex operators $Y(a,z)$ is a Lie superalgebra
$\U_{\ka}(\wh\gl_{m|n})_{\loc}$ with the super-commutator
\beql{commut}
[a_{(r)},b_{(s)}]=\sum_{k\geqslant 0}\binom{r}{k}
\big(a_{(k)}\ts b\big)_{(r+s-k)},
\eeq
which is called the {\it local completion\/}
of the quotient of the universal enveloping algebra
$\U(\wh\gl_{m|n})$ by the ideal generated by $K-\ka$;
see \cite[Sec.~3.5]{fb:va}.

The translation operator is determined by
\beql{transl}
D:1\mapsto 0 \Fand
\big[D, e_{ij}[r]\tss\big]= -r\tss e_{ij}[r-1].
\eeq
The state-field correspondence $Y$
is defined by setting $Y(1,z)=\text{id}$,
\beql{basfi}
Y(e_{ij}[-1],z)=e_{ij}(z):=
\sum_{r\in\ZZ} e_{ij}[r]\tss z^{-r-1},
\eeq
and then extending the map to the whole of $V$
with the use of {\it normal ordering\/}.
Namely, the normally ordered product of homogeneous fields
\ben
a(z)=\sum_{r\in\ZZ}a_{(r)}z^{-r-1}\Fand
b(w)=\sum_{r\in\ZZ}b_{(r)}w^{-r-1}
\een
is the formal power series
\beql{normor}
:a(z)\tss b(w){:}
= a(z)_+\tss b(w)+(-1)^{\deg a\ts\deg b}\ts b(w)\tss a(z)_-,
\eeq
where
\ben
a(z)_+=\sum_{r<0}a_{(r)}z^{-r-1}\Fand
a(z)_-=\sum_{r\geqslant 0}a_{(r)}z^{-r-1}.
\een
This definition extends to an arbitrary number of fields
with the convention that the normal ordering is read
from right to left. Then
\ben
Y(e_{i_1j_1}[-r_1-1]\dots e_{i_mj_m}[-r_m-1],z)
=\frac{1}{r_1!\ts\dots r_m!}
: \di_z^{r_1}e_{i_1j_1}(z)\dots \di_z^{r_m}e_{i_mj_m}(z):.
\een

As defined in the Introduction,
the {\it center\/} of the vertex algebra $V_{n-m}(\gl_{m|n})$
at the critical level $\ka=n-m$ is
\ben
\z(\wh\gl_{m|n})
=\{b\in V_{n-m}(\gl_{m|n})\ |\ \gl_{m|n}[t]\ts b=0\}
\een
which can be identified with a commutative subalgebra
of $\U(t^{-1}\gl_{m|n}[t^{-1}])$.
Elements of $\z(\wh\gl_{m|n})$
are called {\it Segal--Sugawara vectors}.
Due to the commutator formula \eqref{commut},
if $b\in \z(\wh\gl_{m|n})$, then all Fourier
coefficients of the corresponding
field $b(z)=Y(b,z)$ belong to the center
of the Lie superalgebra $\U_{\ka}(\wh\gl_{m|n})_{\loc}$.
These Fourier coefficients are called the {\it Sugawara operators\/}
for $\wh\gl_{m|n}$. In particular, they commute
with the elements of $\wh\gl_{m|n}$ and thus form
a commuting family of $\wh\gl_{m|n}$-endomorphisms
of Verma modules over $\wh\gl_{m|n}$ at the critical level;
cf. \cite{ff:ak}, \cite{gw:ho}, \cite{h:so}.
We will apply the results of Sec.~\ref{sec:ms} to
construct several families of Sugawara operators
for $\wh\gl_{m|n}$.

\subsection{Segal--Sugawara vectors}

Consider the square matrix
\beql{manma}
\tau+\wh E[-1]=\big[\de_{ij}\tau+e_{ij}[-1](-1)^{\bi}\big]
\eeq
with the entries in the universal enveloping algebra
for $\wh\gl_{m|n}\oplus\CC\tau$. The following observation
will play a key role in what follows.

\ble\label{lem:mspm}
The matrix $\tau+\wh E[-1]$ is a Manin matrix.
\ele

\bpf
We have
\ben
\bal
\big[\de_{ij}\tau+e_{ij}[-1](-1)^{\bi},\ts&
\de_{kl}\tau+e_{kl}[-1](-1)^{\bk}\big]\\
{}={}&\de_{ij}\ts e_{kl}[-2](-1)^{\bk}
-\de_{kl}\ts e_{ij}[-2](-1)^{\bi}\\
{}+{}&\de_{kj}\ts e_{il}[-2](-1)^{\bi+\bk}
-\de_{il}\ts e_{kj}[-2](-1)^{(\bi+\bj)(\bk+\bl)+\bi+\bk}.
\eal
\een
This expression remains unchanged after
swapping $i$ and $k$ and multiplying by
$(-1)^{\bi\bj+\bi\bk+\bj\bk}$. Thus,
the matrix elements of $\tau+\wh E[-1]$
satisfy \eqref{srq}.
\epf

\bth\label{thm:traces}
For any $k\geqslant 0$
all coefficients $s_{kl}$ in the expansion
\ben
\str(\tau+\wh E[-1])^k=
s_{k\tss 0}\ts\tau^{k}+s_{k1}\ts\tau^{k-1}
+\dots+s_{kk}
\een
are Segal--Sugawara vectors.
\eth

\bpf
It is sufficient to verify that
for all $i,j$
\beql{annih}
e_{ij}[0]\ts \str(\tau+\wh E[-1])^k=
e_{ij}[1]\ts \str(\tau+\wh E[-1])^k=0
\eeq
in the $\wh\gl_{m|n}$-module
$V_{n-m}(\gl_{m|n})\ot \CC[\tau]$.
We will employ
matrix notation of Sec.~\ref{subsec:mmmt} and consider
the tensor product superalgebra
\ben
\End\CC^{m|n}\ot\dots\ot\End\CC^{m|n}\ot\U
\een
with $k+1$ copies of $\End\CC^{m|n}$ labeled by $0,1,\dots,k$,
where $\U$ stands for the universal enveloping algebra
$\U(\wh\gl_{m|n}\oplus\CC\tau)$.
Set $T=\tau+\wh E[-1]$ and for any integer $r$ introduce
the matrix $\wh E[r]=\big[e_{ij}[r](-1)^{\bi}\big]$. Relations
\eqref{annih} can now be written in the equivalent form
\beql{aanih}
\wh E_0[0]\ts\str_1T_1^{\tss k}=0
\Fand\wh E_0[1]\ts\str_1T_1^{\tss k}=0
\eeq
modulo the left ideal of $\U$ generated by
$\gl_{m|n}[t]$ and $K-n+m$ (or by $K'-1$ in terms of the re-scaled
central element $K'=K/(n-m)$).
In order to verity them, note that by
the commutation relations in the Lie superalgebra
$\wh\gl_{m|n}\oplus\CC\tau$
we have
\begin{align}
[\wh E_0[0],T_1]&=P_{01}T_1-T_1P_{01},
\label{eoan}\\
[\wh E_0[1],T_1]&=\wh E_0[0]+P_{01}\wh E_1[0]-\wh E_1[0]P_{01}
+KP_{01}+\frac{K}{n-m};
\label{eonean}
\end{align}
the last two terms in \eqref{eonean} (and throughout the argument below)
should be re-written in terms of $K'$ in the case $m=n$.
The following identity is well-known:
\beql{compo}
[\wh E_0[0],T_1^{\tss k}]=P_{01}T_1^{\tss k}-T_1^{\tss k}P_{01},
\qquad k=0,1,2,\dots,
\eeq
and it follows immediately from \eqref{eoan}:
\ben
[\wh E_0[0],T_1^{\tss k}]=\sum_{r=1}^k T_1^{\tss r-1}
[\wh E_0[0],T_1] T_1^{\tss k-r}
=\sum_{r=1}^k T_1^{\tss r-1} [P_{01},T_1] T_1^{\tss k-r}
=[P_{01},T_1^{\tss k}].
\een
Now the first relation in \eqref{aanih}
follows by taking the supertrace $\str_1$ on both
sides of \eqref{compo}.

For the proof of the second relation in \eqref{aanih},
suppose first that $m\ne n$ and
use \eqref{eonean} to write
\beql{mapc}
[\wh E_0[1],T_1^{\tss k}]=\sum_{i=1}^k
T_1^{\tss i-1}\Big(\wh E_0[0]+P_{01}\wh E_1[0]-\wh E_1[0]P_{01}
+KP_{01}-\frac{K}{m-n}\Big)\tss T_1^{\tss k-i}.
\eeq
Applying \eqref{compo} and the relation
$\str^{}_{2}T^{k-i}_2P_{12}=T^{k-i}_1$
we can rewrite \eqref{mapc}
modulo the left ideal of $\U$ generated by
$\gl_{m|n}[t]$ as
\ben
\bal[]
[\wh E_0[1],T_1^{\tss k}]&=
\sum_{i=1}^k\Big(
T_1^{\tss i-1}[P_{01},T_1^{\tss k-i}]+KT_1^{\tss i-1}P_{01}
T_1^{\tss k-i}\Big)-\frac{kK}{m-n}\ts T_1^{\tss k-1}\\
{}&+\str_2\sum_{i=1}^k T_1^{\tss i-1}
\big(P_{01}\wh E_1[0]-\wh E_1[0]P_{01}\big)\tss T_2^{\tss k-i}P_{12}.
\eal
\een
Now transform the last summand
using \eqref{compo} to get
\ben
\bal
&\str_2\sum_{i=1}^k T_1^{\tss i-1}
\big(P_{01}\wh E_1[0]-\wh E_1[0]P_{01}\big)\tss T_2^{\tss k-i}P_{12}\\
{}={}&\str_2\sum_{i=1}^k T_1^{\tss i-1}P_{01}[P_{12},T_2^{\tss k-i}]P_{12}
-\str_2\sum_{i=1}^k T_1^{\tss i-1}[P_{12},T_2^{\tss k-i}]P_{01}P_{12}.
\eal
\een
Taking into account the relation
$\str_2 P_{02}=1$, we can simplify this to
\ben
\sum_{i=1}^{k-1}\Big((m-n)T_1^{\tss i-1}P_{01}T_1^{\tss k-i}
-T_1^{\tss i-1}P_{01}\ts\str\ts T^{\tss k-i}+T_1^{\tss i-1}
T_0^{\tss k-i}\Big)-(k-1) T_1^{k-1}.
\een
Combining all the terms and taking the supertrace
$\str_1$ in \eqref{mapc} we derive
\ben
\bal[]
[\wh E_0[1],\str_1\ts T_1^{\tss k}]&=(K+m-n)\Big(T_0^{k-1}
-\frac{k}{m-n}\str\ts T^{k-1}\Big)
+(K+m-n-k+2)\ts T_0^{k-1}\\
{}&+(K+m-n+1)\ts\str_1\sum_{i=2}^{k-1}
P_{01}T_0^{\tss i-1}T_1^{\tss k-i}
-\str_1\sum_{i=2}^{k-1}\ts[T_0^{\tss i-1},T_1^{\tss k-i}].
\eal
\een
Finally, we use Lemmas~\ref{lem:manpr} and \ref{lem:mspm}
to write
\ben
\sum_{i=2}^{k-1}\ts[T_0^{\tss i-1},T_1^{\tss k-i}]
=\sum_{i=2}^{k-1}\ts P_{01}[T_0^{\tss i-1},T_1^{\tss k-i}]
\een
which brings the previous formula to the form
\ben
[\wh E_0[1],\str_1\ts T_1^{\tss k}]=(K+m-n)\Big(2\tss T_0^{k-1}
-\frac{k}{m-n}\str\ts T^{k-1}+\str_1\sum_{i=2}^{k-1}
P_{01}T_0^{\tss i-1}T_1^{\tss k-i}
\Big).
\een
This expression vanishes in
the vacuum module $V_{n-m}(\gl_{m|n})$ as
$K+m-n=0$.

In the case $m=n$ the same calculation applies with the use
of \eqref{eonean} rewritten in terms of $K'=K/(n-m)$. The critical value
of $K'$ is understood as being equal to $1$.
\epf

The following corollary is immediate
from Theorems~\ref{thm:ident} and \ref{thm:traces}.

\bco\label{cor:ss}
All the coefficients
$\si_{k\tss l}, h_{k\tss l},
b_{k\tss l}\in\U(t^{-1}\gl_{m|n}[t^{-1}])$
in the expansions
\ben
\bal
\str^{}_{1,\dots,k}\ts A_k T_1\dots T_k&=
\si_{k\tss 0}\tau^k+\si_{k1}\tau^{k-1}+\dots+\si_{kk},\\
\str^{}_{1,\dots,k}\ts H_k T_1\dots T_k&=
h_{k\tss 0}\tau^k+h_{k1}\tau^{k-1}+\dots+h_{kk},\\
\Ber\big(1+u\tss T\big)&=
\sum_{k=0}^{\infty} \sum_{l=0}^k b_{k\tss l}\ts u^k\tau^{k-l}
\eal
\een
are Segal--Sugawara vectors. Moreover, $b_{k\tss l}=\si_{k\tss l}$
for all $k$ and $l$.
\qed
\eco

\bre\label{rem:comp}
(i)\quad In the case $n=0$ the formulas of
Theorem~\ref{thm:traces} and Corollary~\ref{cor:ss}
reproduce the explicit constructions
of Segal--Sugawara vectors for the affine Kac--Moody algebra $\wh\gl_m$;
see \cite{cm:ho} and \cite{ct:qs}. The Berezinian $\Ber\big(1+u\tss T\big)$
turns into the column determinant $\cdet\tss\big(1+u\tss T\big)$
which is a polynomial in $u$, and the proof given
in \cite{cm:ho} deals with its leading coefficient
$\cdet\ts T$. The use of the matrix techniques
in the above proof of Theorem~\ref{thm:traces} thus provided a simpler proof
of those results for $\wh\gl_m$.

(ii)\quad
Denote by $\bar s_{k\tss k}$, $\bar h_{k\tss k}$ and
$\bar b_{k\tss k}$ the respective images
of the elements $s_{k\tss k}$, $h_{k\tss k}$ and
$b_{k\tss k}$
in the associated graded algebra
$\text{gr}\ts\U(t^{-1}\gl_{m|n}[t^{-1}])
\cong \Sr(t^{-1}\gl_{m|n}[t^{-1}])$.
Each of the families $\bar s_{k\tss k}$, $\bar h_{k\tss k}$ and
$\bar b_{k\tss k}$ with $k\geqslant 1$ is the image
of a generating set of the algebra
of invariants $\Sr(\gl_{m|n})^{\gl_{m|n}}$ under the embedding
$\Sr(\gl_{m|n})\hookrightarrow
\Sr(t^{-1}\gl_{m|n}[t^{-1}])$ defined by
the assignment $X\mapsto X[-1]$.
In the case $n=0$ the families of Segal--Sugawara vectors
with such a property are known as
{\it complete sets of Segal--Sugawara vectors\/}
(see \cite{cm:ho}, \cite{ff:ak}, \cite{gw:ho}, \cite{h:so})
so that this terminology can be extended to the super case.
It is natural to suppose that each of the families
$D^r s_{k\tss k}$, $D^r h_{k\tss k}$ and
$D^r b_{k\tss k}$ with $r\geqslant 0$
and $k\geqslant 1$ generates the center $\z(\wh\gl_{m|n})$
of the vertex algebra $V_{n-m}(\gl_{m|n})$; cf.
\cite{f:wm}. However, we do not have a proof of this
conjecture. If both $m$ and $n$ are positive integers, then
none of the families is algebraically independent;
cf. Theorem~\ref{thm:gensv} below.
\qed
\ere

The application of the state-field correspondence map $Y$
to the Segal--Sugawara vectors produces elements of the
center of the local completion $\U_{n-m}(\wh\gl_{m|n})_{\loc}$
at the critical level $\ka=n-m$. Hence, Theorem~\ref{thm:traces}
and Corollary~\ref{cor:ss} provide explicit formulas
for the corresponding Sugawara operators for $\wh\gl_{m|n}$.

Recall the matrix
$\wh E(z)=[e_{ij}(z)(-1)^{\bi}]$, where the fields
$e_{ij}(z)$ are defined in \eqref{basfi}.
Set $T(z)=\di_z+\wh E(z)$, where $\di_z$ is understood
as a scalar matrix of size $m+n$.

\bco\label{cor:sop}
All Fourier coefficients of the fields
$s_{kl}(z)$, $b_{kl}(z)$, $\si_{kl}(z)$ and $h_{kl}(z)$
defined by the decompositions
\ben
\bal
:\str\ts T(z)^k:{}&
= s_{k\tss 0}(z)\ts\di_z^{\tss k}+s_{k1}(z)\ts\di_z^{\tss k-1} +\dots
+s_{kk}(z),\\
:\Ber\big(1+u\ts T(z)\big):{}&=
\sum_{k=0}^{\infty} \sum_{l=0}^k b_{kl}(z)\ts u^k\di_z^{\ts k-l},\\
:\str^{}_{1,\dots,k}\ts A_k T_1(z)\dots T_k(z):{}&=
\si_{k\tss 0}(z)\ts\di_z^{\tss k}+\si_{k1}(z)\ts\di_z^{\tss k-1}
+\dots+\si_{kk}(z),\\
:\str^{}_{1,\dots,k}\ts H_k T_1(z)\dots T_k(z):{}&=
h_{k\tss 0}(z)\ts\di_z^{\tss k}+h_{k1}(z)\ts\di_z^{\tss k-1}
+\dots+h_{kk}(z)
\eal
\een
are Sugawara operators for $\wh\gl_{m|n}$.
Moreover, $\si  _{kl}(z)=b_{kl}(z)$ for all $k$ and $l$.
\qed
\eco

\subsection{Commutative subalgebras and super Gaudin
Hamiltonians}
\label{subsec:cssg}

By the vacuum axiom of a vertex algebra, the application
of any of the fields
$s_{kl}(z)$, $b_{kl}(z)$, $\si_{kl}(z)$ and $h_{kl}(z)$
introduced in Corollary~\ref{cor:sop}
to the vacuum vector yields formal power series in $z$
with coefficients in the universal enveloping algebra
$\U(t^{-1}\gl_{m|n}[t^{-1}])$. Moreover, since
the Fourier coefficients of the fields belong to the center
of the local completion $\U_{n-m}(\wh\gl_{m|n})_{\loc}$,
all coefficients of the power series will belong
to the center $\z(\wh\gl_{m|n})$ of the vertex algebra
$V_{n-m}(\gl_{m|n})$. In particular, these coefficients
generate a commutative subalgebra
of $\U(t^{-1}\gl_{m|n}[t^{-1}])$. Explicitly, they can be given
by the same expansions as in Corollary~\ref{cor:sop}
by omitting the normal ordering signs and by
replacing $T(z)$ with the matrix
$T(z)_+=\di_z+\wh E(z)_+$ with
$\wh E(z)_+=[(-1)^{\bi}e_{ij}(z)_+]$ and
\beql{serpos}
e_{ij}(z)_+=\sum_{r=0}^{\infty} e_{ij}[-r-1]\tss z^r.
\eeq

Now we produce the corresponding families of commuting
elements of $\U(\gl_{m|n}[t])$.
The commutation relations of the Lie superalgebra
$t^{-1}\gl_{m|n}[t^{-1}]$ can be written in terms of the
series \eqref{serpos} as
\begin{align}\non
(z-w)\ts[e_{ij}(z)_+,e_{kl}(w)_+]
&=\de_{kj}\big(e_{il}(z)_+-e_{il}(w)_+\big)\\
\label{relpl}
{}&-\de_{il}\big(e_{kj}(z)_+-e_{kj}(w)_+\big)
(-1)^{(\bi+\bj)(\bk+\bl)}.
\end{align}
Setting
\ben
e_{ij}(z)_-=\sum_{r=0}^{\infty} e_{ij}[r]\tss z^{-r-1},
\een
we find that
\begin{align}\non
(z-w)\ts[e_{ij}(z)_-,e_{kl}(w)_-]
{}={}&-{}\de_{kj}\big(e_{il}(z)_--e_{il}(w)_-\big)\\
\label{relmn}
&+\de_{il}\big(e_{kj}(z)_--e_{kj}(w)_-\big)
(-1)^{(\bi+\bj)(\bk+\bl)}.
\end{align}

The commuting families of
elements of $\U(t^{-1}\gl_{m|n}[t^{-1}])$
obtained by the application of Corollary~\ref{cor:sop}
are expressed as coefficients of certain differential
polynomials in the series $e_{ij}(z)_+$
(i.e., polynomials in $\di_z^{\tss k}e_{ij}(z)_+$).
Therefore, comparing \eqref{relpl} and \eqref{relmn},
we can conclude that the same differential polynomials
with the $e_{ij}(z)_+$ respectively replaced by
$-e_{ij}(z)_-$, generate a commutative subalgebra of
$\U(\gl_{m|n}[t])$. Indeed,
suppose that $A(z)$ is a differential
polynomial in the $e_{ij}(z)_+$ and $B(w)$ is
a differential polynomial in the $e_{ij}(w)_+$.
A relation of the form
$[A(z),B(w)]=0$ is a consequence of the commutation
relations between the series \eqref{relpl} where the
actual expansions of the $e_{ij}(z)_+$
and $e_{ij}(w)_+$ as power series
in $z$ and $w$ need not be used.
This relies on the easily verified property that given
any total ordering on the set of series $\di_z^{\tss k}e_{ij}(z)_+$,
the corresponding ordered monomials
are linearly independent over the polynomial ring $\CC[z]$.
Thus we arrive at the following
corollary, where we use the notation
\ben
L(z)=\di_z-\wh E(z)_-,\qquad
\wh E(z)_-=[(-1)^{\bi}e_{ij}(z)_-].
\een

\bco\label{cor:commsa}
All coefficients of the series
$S_{kl}(z)$, $B_{kl}(z)$, $\Sigma_{kl}(z)$ and $H_{kl}(z)$
defined by the decompositions
\ben
\bal
\str\ts L(z)^k{}&
= S_{k\tss 0}(z)\ts\di_z^{\tss k}+
S_{k1}(z)\ts\di_z^{\tss k-1} +\dots
+S_{kk}(z),\\
\Ber\big(1+u\ts L(z)\big){}&=
\sum_{k=0}^{\infty} \sum_{l=0}^k B_{kl}(z)\ts u^k\di_z^{\ts k-l},\\
\str^{}_{1,\dots,k}\ts A_k L_1(z)\dots L_k(z){}&=
\Sigma_{k\tss 0}(z)\ts\di_z^{\tss k}+
\Sigma_{k1}(z)\ts\di_z^{\tss k-1}
+\dots+\Sigma_{kk}(z),\\
\str^{}_{1,\dots,k}\ts H_k L_1(z)\dots L_k(z){}&=
H_{k\tss 0}(z)\ts\di_z^{\tss k}+H_{k1}(z)\ts\di_z^{\tss k-1}
+\dots+H_{kk}(z)
\eal
\een
generate a commutative subalgebra of $\U(\gl_{m|n}[t])$.
Moreover, $\Sigma_{kl}(z)=B_{kl}(z)$ for all
values of $k$ and $l$.
\qed
\eco

Note also that given any set of complex parameters
$\la_1,\dots,\la_{m+n}$ the commutation relations
\eqref{relmn} remain valid after the replacement
\beql{replala}
e_{ij}(z)_-\mapsto \de_{ij}\la_i+e_{ij}(z)_-
\eeq
Thus, we obtain more general families
of commutative subalgebras of $\U(\gl_{m|n}[t])$.

\bco\label{cor:repla}
Given a set of complex parameters
$\la_1,\dots,\la_{m+n}$,
all coefficients of the series
$S_{kl}(z)$, $B_{kl}(z)$, $\Sigma_{kl}(z)$ and $H_{kl}(z)$
defined by the decompositions
of Corollary~\ref{cor:commsa} with the series
$e_{ij}(z)_-$ replaced as in \eqref{replala}
generate a commutative subalgebra of $\U(\gl_{m|n}[t])$.
\qed
\eco

\bre\label{rem:commsab}
(i)\quad A more general replacement
$e_{ij}(z)_-\mapsto K_{ij}+e_{ij}(z)_-$
with appropriately defined elements $K_{ij}$
preserves the commutation relations \eqref{relmn} as well,
thus leading to even more general family of commutative
subalgebras of $\U(\gl_{m|n}[t])$.

(ii)\quad In the particular case $n=0$ the elements
of the commutative subalgebras of $\U(\gl_m[t])$
from Corollaries~\ref{cor:commsa} and \ref{cor:repla}
were originally constructed in \cite{t:qg};
see also \cite{cm:ho}, \cite{ct:qs}, \cite{mtv:be},
\cite{mtv:ba}. In particular, the Berezinian turns into
a determinant, and generators of the commutative
subalgebra are found by the decomposition of the column determinant
\ben
\cdet\left[\begin{matrix}
\di_z-e_{11}(z)_-&-e_{12}(z)_-&\dots&-e_{1m}(z)_-\\
-e_{21}(z)_-&\di_z-e_{22}(z)_-&\dots&-e_{2m}(z)_-\\
\vdots&\vdots& \ddots&\vdots     \\
-e_{m1}(z)_-&-e_{m2}(z)_-&\dots&\di_z-e_{mm}(z)_-
                \end{matrix}\right].
\een
\vskip-1.3\baselineskip
\qed
\ere

The commutative subalgebras of $\U(\gl_{m|n}[t])$
can be used to construct higher Gaudin Hamiltonians
following the same scheme as in the even case; see
\cite{ffr:gm}, \cite{mtv:be}, \cite{t:qg}.
More precisely, given a finite-dimensional $\gl_{m|n}$-module
$M$ and a complex number $a$ we can define
the corresponding evaluation $\gl_{m|n}[t]$-module $M_a$.
As a vector superspace, $M_a$ coincides with $M$,
while the action of the elements of the Lie
superalgebra is given by $e_{ij}[r]\mapsto e_{ij}\ts a^r$ for
$r\geqslant 0$, or equivalently,
\ben
e_{ij}(z)_-\mapsto \frac{e_{ij}}{z-a}.
\een
Now consider certain
finite-dimensional $\gl_{m|n}$-modules $M^{(1)},\dots,M^{(k)}$
and let $a_1,\dots,a_k$ be complex parameters.
Then the tensor product of the evaluation modules
\beql{gaudmo}
M^{(1)}_{a_1}\ot\dots\ot M^{(k)}_{a_k}
\eeq
becomes a $\gl_{m|n}[t]$-module via the iterated
comultiplication map \eqref{comult}.
Then the images of the matrix elements of the matrix
$L(z)$ are found by
\ben
\ell_{ij}(z)=\de_{ij}\di_z
-(-1)^{\bi}\sum_{r=1}^k\frac{e_{ij}^{(r)}}{z-a_r},
\een
where $e_{ij}^{(r)}$ denotes the image of $e_{ij}$
in the $\gl_{m|n}$-module $M^{(r)}$. So replacing
$L(z)$ by the matrix $\Lc(z)=[\ell_{ij}(z)]$
in the formulas of Corollary~\ref{cor:commsa}
we obtain a family of commuting operators
in the module \eqref{gaudmo}, thus producing higher
Gaudin Hamiltonians associated with $\gl_{m|n}$.
In particular, such families are provided by
the coefficients of the series
defined by the expansion of the Berezinian
$\Ber\big(1+u\ts \Lc(z)\big)$
and the supertrace
\ben
\str\ts \Lc(z)^k
=\Sc_{k\tss 0}(z)\ts\di_z^{\tss k}+
\Sc_{k1}(z)\ts\di_z^{\tss k-1} +\dots
+\Sc_{kk}(z).
\een
The quadratic Gaudin Hamiltonian
$\Hc(z)=\Sc_{22}(z)$ can be written explicitly as
\ben
\Hc(z)=\sum_{r,s=1}^k\frac{1}{(z-a_r)(z-a_s)}
\sum_{i,j}e_{ij}^{(r)}e_{ji}^{(s)}(-1)^{\bj}+
\sum_{r=1}^k\frac{1}{(z-a_r)^2}
\sum_{i}e_{ii}^{(r)}.
\een
Assuming further that the parameters $a_i$ are all distinct
and setting
\ben
\Hc^{(r)}=\sum_{s\ne r}\frac{1}{a_r-a_s}
\sum_{i,j}e_{ij}^{(r)}e_{ji}^{(s)}(-1)^{\bj}
\een
we can also write the Hamiltonian as
\ben
\Hc(z)=2\tss\sum_{r=1}^k \frac{\Hc^{(r)}}{z-a_r}
+\sum_{r=1}^k\frac{\Delta^{(r)}}{(z-a_r)^2},
\een
where $\Delta^{(r)}$ denotes the eigenvalue of the Casimir
element $\sum e_{ij}e_{ji}(-1)^{\bj}+\sum e_{ii}$
of $\gl_{m|n}$
in the representation $M^{(r)}$; cf. \cite{ffr:gm}, \cite{km:co}.

More general families of commuting elements in $\U(\gl_{m|n}[t])$
and higher Gaudin Hamiltonians
can be constructed by using extra parameters $\la_i$
or $K_{ij}$ as in Corollary~\ref{cor:repla}
and Remark~\ref{rem:commsab}(i).

\subsection{Singular vectors in Verma modules}
\label{subsec:svvm}

We now recall a result from \cite{g:gv}
and formulate it for the affine Lie superalgebra
of type $A(m-1\tss|\tss n-1)$.
Consider the Lie superalgebra
$\wh\g=\wh\gl_{m|n}\oplus \CC d$, so that
\ben
\wh\g=\gl_{m|n}[t,t^{-1}]\oplus\CC K\oplus \CC d,
\een
with the commutation relations \eqref{commrel},
where the elements $K$ and $d$ are even, $K$ is central and
\beql{dco}
\big[d,e_{ij}[r]\tss\big]=r\ts e_{ij}[r].
\eeq
Denote by $\n_{-},\h$ and $\n$ the subalgebras of $\gl_{m|n}$
spanned by the lower-triangular, diagonal and upper triangular
matrices, respectively. We have the triangular decomposition
\ben
\wh\g=\wh\n_{-}\oplus \wh\h\oplus \wh\n,
\een
where $\wh\h$ is the Cartan subalgebra of $\wh\g$ defined by
\ben
\wh\h=\h\oplus\CC K\oplus \CC d,
\een
while
\ben
\wh\n_{-}=\n_{-}[t^{-1}]\oplus t^{-1}\h[t^{-1}]
\oplus t^{-1}\n[t^{-1}]
\Fand
\wh\n=t\ts\n_{-}[t]\oplus t\ts\h[t]
\oplus \n[t].
\een
We have the direct sum decomposition
\ben
\U(\wh\n_{-})=\U(t^{-1}\h[t^{-1}])\oplus
\Big(\U(\wh\n_{-})\n_{-}[t^{-1}]+t^{-1}\n[t^{-1}]\U(\wh\n_{-})\Big).
\een
Introduce the projection to the first component in the direct sum
\beql{hc}
\HC:\U(\wh\n_{-})\to\U(t^{-1}\h[t^{-1}])
\eeq
so that the kernel of $\HC$ coincides with the second component.

Write the dual vector space $\wh\h^*$ to the Cartan subalgebra
in the form
\ben
\wh\h^*=\h^*\oplus\CC \La_0\oplus \CC \de,
\een
where
\ben
\La_0(\h\oplus\CC d)=0,\qquad \de(\h\oplus\CC K)=0,
\qquad \La_0(K)=\de(d)=1
\een
and $\nu(K)=\nu(d)=0$ for all $\nu\in \h^*$.
Let $\ve_1,\dots,\ve_{m+n}\in \h^*$ denote the basis of
$\h^*$ dual to the standard basis of $\h$ so that
$\ve_i(e_{jj})=\de_{ij}$.

Given an element
$\la\in\wh\h^*$,
the {\it Verma module\/} $M(\la)$ with the highest weight $\la$
is defined as the induced module
\ben
M(\la)=\U(\wh\g)\ot_{\U(\wh\h\oplus \wh\n)}\CC_{\la},
\een
where $\CC_{\la}$ is the one-dimensional $(\wh\h\oplus \wh\n)$-module
with the basis vector $1_{\la}$ such that
\ben
\wh\n\ts 1_{\la}=0\Fand h \ts 1_{\la}=\la(h)\ts 1_{\la},\quad
h\in\wh\h.
\een
We will identify $M(\la)$
with $\U(\wh\n_{-})$ as a vector space.

The highest weight $\la$ is {\it at the critical level\/},
if $\la(K)=n-m$.
The vector space $\wh\h^*$
is equipped with the bilinear form $(\ ,\ )$ defined by
\begin{multline}
(\nu+a\La_0+b\tss\de,\mu+a'\La_0+b'\de)\\
=\nu_1\mu_1+\dots+\nu_m\mu_m-\nu_{m+1}\mu_{m+1}
-\dots-\nu_{m+n}\mu_{m+n}
+a\tss b'+b\tss a',
\non
\end{multline}
where $\nu=\nu_1\ve_1+\dots+\nu_{m+n}\ve_{m+n}$
and $\mu=\mu_1\ve_1+\dots+\mu_{m+n}\ve_{m+n}$.
The condition that $\la$ is at the critical level
can be written as $(\la+\rho,\de)=0$, where
the element $\rho\in\wh\h^*$ is defined by
\ben
\rho=(m-1)\tss\ve_1+\dots+\ve_{m-1}-\ve_{m+2}-\dots-(n-1)\tss\ve_{m+n}
+(m-n)\La_0.
\een
The root system $\De\subset \wh\h^*$ of $\wh\g$ is formed by the
elements $k\tss\de+\ve_i-\ve_j$ with $k\in\ZZ$ and $i\ne j$ and by
the elements $k\tss\de$ with $k\in\ZZ$, $k\ne 0$.
The real positive roots are given by
\ben
\De^+_{re}=\{k\tss\de+\ve_i-\ve_j\ |\ k>0,\quad i\ne j\}\cup
\{\ve_i-\ve_j\ |\ i< j\}.
\een
An element $\la\in\wh\h^*$ is called a {\it generic
critical weight\/} if $\la$ is at the critical level
and
\beql{gener}
(\la+\rho,\al)\ne \frac{p}{2}\ts(\al,\al)\qquad\text{for all}\quad
p=1,2,\dots\fand \al\in\De^+_{re}.
\eeq

A {\it singular vector\/} of the Verma module
is an element $v\in M(\la)$ such that
$\wh\n\ts v=0$. If $\la$ is a generic
critical weight, then
the subspace of singular vectors $M(\la)^{\wh\n}$ has
the weight space decomposition of the form
\ben
M(\la)^{\wh\n}=\bigoplus_{k\geqslant 0} M(\la)^{\wh\n}_{\la-k\tss\de}.
\een
Moreover, $M(\la)^{\wh\n}$ possesses
an algebra structure which is described in
\cite[Theorem~1.1]{g:gv}:
the restriction of the projection $\HC$ defined in \eqref{hc}
to $M(\la)^{\wh\n}$ gives an algebra isomorphism
\beql{isomsv}
\HC:M(\la)^{\wh\n}\cong\U(t^{-1}\h[t^{-1}]).
\eeq
This result was applied in \cite{g:gv} to
obtain the character formula for the simple quotient
of $M(\la)$
thus proving the Kac--Kazhdan conjecture; see also
an earlier paper \cite{ik:wm} for another proof
of this conjecture for affine Lie superalgebras of type
$A(m-1\ts|\ts n-1)$.

The next lemma will be needed for the proof
of Theorem~\ref{thm:gensv} below.
For this lemma
we regard $\la_1,\dots,\la_{m+n}$ as independent variables.
Consider the square matrix $\La=[\La_{ki}]$
with $1\leqslant k,i\leqslant m+n$
whose entries are
polynomials in the $\la_i$ and a complex parameter
$r$ defined by the formulas
\beql{matele}
\bal
\La_{ki}=\sum_{a+b+c=k-1}
&h_a(\la_{m+1}+r+k-a,\dots,\la_{m+n}+r+k-a)\\
{}&\times{}e_b(\la_{i+1}-r-c-1,\dots,\la_m-r-c-1)\ts
e_c(\la_1,\dots,\la_{i-1})
\eal
\eeq
for $i=1,\dots,m$, and
\beql{matere}
\bal
\La_{ki}=\sum_{a+b+c=k-1}
&h_a(\la_i+r+k-a,\dots,\la_{m+n}+r+k-a)\\
{}&\times{}h_b(\la_{m+1}+c,\dots,\la_i+c)\ts
e_c(\la_1,\dots,\la_m)
\eal
\eeq
for $i=m+1,\dots,m+n$ with the summation taken over nonnegative
integers $a,b,c$,
and we used the following
shifted versions of the complete and elementary
symmetric polynomials:
\ben
h_a(x_1,\dots,x_l)
=\sum_{l\geqslant i_1\geqslant\dots\geqslant i_a\geqslant 1}
(x_{i_1}+a-1)\dots (x_{i_{a-1}}+1)\ts x_{i_a}
\een
and
\ben
e_a(x_1,\dots,x_l)
=\sum_{l> i_1>\dots> i_a> 1}
(x_{i_1}-a+1)\dots (x_{i_{a-1}}-1)\ts x_{i_a}.
\een

\ble\label{lem:det}
We have
\ben
\bal
\det\La=\prod_{1\leqslant i<j\leqslant m}
(\la_i-\la_j+j-i+r)&\prod_{m+1\leqslant i<j\leqslant m+n}
(\la_j-\la_i-j+i-r)
\\
{}\times{}&\prod_{1\leqslant i\leqslant m<j\leqslant m+n}
(\la_i+\la_j-i-j+2m+1).
\eal
\een
\ele

\bpf
Note that all entries in the top row of $\La$ are equal to
$1$. Then replacing the $i$-th column $col_i$ of $\La$
by the difference
$col_{i+1}-col_i$ for $i=1,\dots,m+n-1$
we find that $\det\La$ equals the determinant of another
square matrix of size $m+n-1$. All elements
in each column of this matrix
has a common linear factor. By taking out all these factors
we are left with a matrix whose entries
in the top row are equal to $1$ and we can repeat the same
procedure.

More precisely,
define a sequence of matrices $\La^{(s)}=[\La^{(s)}_{ki}]$
with $s=0,1,\dots,m+n-1$ so that $\La^{(s)}$ is a square
matrix of size $m+n-s$. We set $\La^{(0)}=\La$ and
define the entries of $\La^{(s)}$ inductively.
For $s\geqslant 1$ we set
\beql{reclem}
\La^{(s)}_{ki}=\frac{\La^{(s-1)}_{k+1,i+1}-\La^{(s-1)}_{k+1,i}}
{\la_{i}-\la_{i+s}+s+r}
\eeq
if $1\leqslant i\leqslant m-s$,
\beql{recmid}
\La^{(s)}_{ki}=\frac{\La^{(s-1)}_{k+1,i+1}-\La^{(s-1)}_{k+1,i}}
{\la_{i}+\la_{i+s}-2i-s+2m+1}
\eeq
if $m-s+1\leqslant i\leqslant m$, and
\beql{recrem}
\La^{(s)}_{ki}=\frac{\La^{(s-1)}_{k+1,i+1}-\La^{(s-1)}_{k+1,i}}
{\la_{i+s}-\la_{i}-s-r}
\eeq
if $m+1\leqslant i\leqslant m+n-s$. It is straightforward
to verify that these relations are satisfied by
the matrix elements
\ben
\bal
\La^{(s)}_{ki}=\sum_{a+b+c=k-1}
&h_a(\la_{m+1}+r+k+s-a,\dots,\la_{m+n}+r+k+s-a)\\
{}&\times{}e^{}_b(\la_{i+s+1}-r-c-s-1,\dots,\la_m-r-c-s-1)\ts
e^{}_c(\la_1,\dots,\la_{i-1})
\eal
\een
for $1\leqslant i\leqslant m-s$,
\ben
\bal
\La^{(s)}_{ki}=\sum_{a+b+c=k-1}
&h_a(\la_{m+1}+r+k+s-a,\dots,\la_{m+n}+r+k+s-a)\\
{}&\times{}h^{}_b(\la_{m+1}+c+m-i+1,\dots,\la_{i+s}+c+m-i+1)\ts
e^{}_c(\la_1,\dots,\la_{i-1})
\eal
\een
for $m-s+1\leqslant i\leqslant m$ and
\ben
\bal
\La^{(s)}_{ki}=\sum_{a+b+c=k-1}
&h_a(\la_i+r+k+s-a,\dots,\la_{m+n}+r+k+s-a)\\
{}&\times{}h^{}_b(\la_{m+1}+c,\dots,\la_{i+s}+c)\ts
e^{}_c(\la_1,\dots,\la_m)
\eal
\een
for $m+1\leqslant i\leqslant m+n-s$. Indeed,
if $1\leqslant i\leqslant m-s$,
then the difference $\La^{(s-1)}_{k+1,i+1}-\La^{(s-1)}_{k+1,i}$
can be found by calculating
\begin{multline}
\non
\sum_{b+c=d}
e^{}_b(\la_{i+s+1}-r-c-s,\dots,\la_m-r-c-s)\ts
e^{}_c(\la_1,\dots,\la_{i})\\
{}-\sum_{b+c=d}
e^{}_b(\la_{i+s}-r-c-s,\dots,\la_m-r-c-s)\ts
e^{}_c(\la_1,\dots,\la_{i-1})
\end{multline}
with $d=k-a-1$. Using the recurrence relations
\ben
e_a(x_1,\dots,x_l)=e_a(x_2,\dots,x_l)
+e_{a-1}(x_2-1,\dots,x_l-1)\ts x_1
\een
and
\ben
e_a(x_1,\dots,x_l)=e_a(x_1,\dots,x_{l-1})
+(x_l-a+1)\ts e_{a-1}(x_1,\dots,x_{l-1})
\een
we bring the difference of the sums to the form
\ben
(\la_{i}-\la_{i+s}+s+r)
\sum_{b+c=d-1}
e^{}_b(\la_{i+s+1}-r-c-s-1,\dots,\la_m-r-c-s-1)\ts
e^{}_c(\la_1,\dots,\la_{i-1})
\een
thus verifying \eqref{reclem}. Checking that the $\La^{(s)}_{ki}$
satisfy \eqref{recmid} and \eqref{recrem} is quite similar
by using analogous recurrence relations
for the polynomials $h_a(x_1,\dots,x_l)$ together with
the relation
\ben
h_a(x_1,\dots,x_l)=h_a(x_1+1,\dots,x_{l-1}+1)
+(x_l-l+1)\ts h_{a-1}(x_1+1,\dots,x_{l}+1)
\een
which follows from the
well-known formula for the generating series of these
polynomials:
\ben
1+\sum_{a=1}^{\infty}\frac{h_a(x_1,\dots,x_l)}
{(u+1)\dots (u+a)}=\frac{u(u-1)\dots (u-l+1)}
{(u-x_1-l+1)\dots (u-x_l)}.
\een

Since the determinant of the $1\times 1$ matrix $\La^{(m+n-1)}=[1]$
is equal to $1$, taking the product of the linear factors
occurring in the denominators in \eqref{reclem},
\eqref{recmid} and \eqref{recrem} for all values
of $s$, we get the desired
formula for $\det\La$.
\epf

Now recall that the Fourier coefficients of all fields
defined in Corollary~\ref{cor:sop} belong
to the center of the
local completion $\U_{n-m}(\wh\gl_{m|n})_{\loc}$
of the universal enveloping algebra
$\U(\wh\gl_{m|n})$. Each of these coefficients
is a well-defined operator
in the Verma module $M(\la)$ at the critical level,
preserving the subspace $M(\la)^{\wh\n}$
of singular vectors. We now aim to produce
algebraically independent generators of
$M(\la)^{\wh\n}$ regarded as a commutative algebra
for generic critical weights $\la$.

Recall the fields $\si_{kl}(z)$ defined
in Corollary~\ref{cor:sop} and introduce their
Fourier coefficients by the expansions
\ben
\si_{kl}(z)=\sum_{r\in\ZZ}\si_{kl}[r]z^{-r-l}.
\een
Then $\si_{kl}[r]\ts 1_{\la}$ is a singular
vector of $M(\la)$ of weight $\la+r\de$
(this vector is zero for $r>0$).

\bth\label{thm:gensv}
Suppose that $\la$ is a generic critical weight.
Then
\ben
M(\la)^{\wh\n}
=\CC[\si_{kk}[r]\ |\  k=1,\dots,m+n,\  r<0\tss]\ts 1_{\la}
\een
so that any singular vector of $M(\la)$ can be expressed
as a unique polynomial in
the operators
$\si_{kk}[r]$ with $k=1,\dots,m+n$ and $r<0$
applied to the highest vector.
\eth

\bpf
Denote by $\si^{\tss\la}_{kl}[r]$
the image of the element $\si_{kl}[r]\ts 1_{\la}$
under the isomorphism $\HC$ defined in \eqref{isomsv}.
It will be sufficient
to show that the elements $\si^{\tss\la}_{kk}[r]$
with $k=1,\dots,m+n$ and $r<0$
are algebraically
independent generators of the algebra $\U(t^{-1}\h[t^{-1}])$.
We will use the explicit formulas for the operators
$\si_{kk}[r]$ implied by Proposition~\ref{prop:extsh}.
Keeping the same notation and
applying the first relation of that proposition we get
\beql{srte}
\bal
:\str^{}_{1,\dots,k}\ts A_k T_1(z)\dots T_k(z):
&=\sum_{I}\frac{1}{\al_{m+1}!\dots\al_{m+n}!}
\sum_{\si\in\Sym_k}\sgn\si\cdot\vp(\si,\si I,I)\\
{}&\times{}:\big(\di_z+\wh E(z)\big)_{i_{\si(1)}i_1}\dots
\big(\di_z+\wh E(z)\big)_{i_{\si(k)}i_k}:(-1)^{\ga(\si I,I)}.
\eal
\eeq
Given $I$ and $\si$, use the definition
of the normal ordering \eqref{normor}
and apply the corresponding summand
to the highest vector:
\begin{multline}
\non
:\big(\di_z+\wh E(z)\big)_{i_{\si(1)}i_1} C(z):\ts 1_{\la}\\
=\big(\de_{i_{\si(1)}i_1}\di_z
+(-1)^{\bi_{\si(1)}}\ts
e_{i_{\si(1)}i_1}(z)_+\big)\ts C(z)\ts 1_{\la}
+(-1)^{\bi_1}\ts
C(z)\ts e_{i_{\si(1)}i_1}(z)_-\ts 1_{\la},
\end{multline}
where we set
\ben
C(z)={}:\big(\di_z+\wh E(z)\big)_{i_{\si(2)}i_2}\dots
\big(\di_z+\wh E(z)\big)_{i_{\si(k)}i_k}:.
\een
We have
\beql{applhw}
e_{i_{\si(1)}i_1}(z)_-\ts 1_{\la}=
\sum_{r\geqslant 0} z^{-r-1}\ts e_{i_{\si(1)}i_1}[r] \ts 1_{\la}.
\eeq
Since $i_1\geqslant\dots\geqslant i_l\geqslant m+1>i_{l+1}>\dots>i_k$,
this expression is zero unless $i_{\si(1)}=i_1$. In the latter
case it equals $z^{-1}\ts \la_{i_1}\ts 1_{\la}$, where the $\la_i$
are the components of the highest weight
$\la=\la_1\ve_1+\dots+\la_{m+n}\ve_{m+n}+(n-m)\La_0+b\ts\de$
(the value of $b$ is irrelevant).
On the other hand,
\ben
e_{i_{\si(1)}i_1}(z)^{}_+\ts C(z)\ts 1_{\la}
=\sum_{r< 0} z^{-r-1}\ts e_{i_{\si(1)}i_1}[r]\ts C(z) \ts 1_{\la}.
\een
The coefficient of each power of $z$ belongs to
the kernel of the projection \eqref{hc} unless
$i_{\si(1)}=i_1$. Proceeding by induction and using
the same argument for $C(z)\tss 1_{\la}$ we conclude that
after the application to the highest weight vector $1_{\la}$,
the summation in \eqref{srte} can be restricted to those
permutations $\si$ which stabilize the multiset $I$
so that $\si\ts I=I$. The number of such permutations
is $\al_{m+1}!\dots\al_{m+n}!$. Furthermore, by \eqref{phinot}
for any such $\si$ the sign $\vp(\si,I,I)$ coincides
with $\sgn\si$ so that under
the isomorphism
$\HC$ defined in \eqref{isomsv} we have
\begin{multline}
\label{proddz}
:\str^{}_{1,\dots,k}\ts A_k T_1(z)\dots T_k(z):{}1_{\la}\mapsto\\
\sum_I \big({-}\di_z+\la_{i_1}z^{-1}+e_{i_1i_1}(z)_+\big)\dots
\big({-}\di_z+\la_{i_l}z^{-1}+e_{i_li_l}(z)_+\big)\\
{}\times{}\big(\ts\di_z+\la_{i_{\tss l+1}}z^{-1}
+e_{i_{\tss l+1}i_{\tss l+1}}(z)_+\big)
\dots
\big(\ts\di_z+\la_{i_k}z^{-1}+e_{i_ki_k}(z)_+\big),
\end{multline}
summed over multisets
$I=\{i_1\geqslant\dots\geqslant i_l\geqslant m+1>i_{l+1}>\dots>i_k\}$
with $l=0,\dots,k$. Thus the elements
$\si^{\tss\la}_{kl}[r]\in\U(t^{-1}\h[t^{-1}])$ are found by writing
the differential operator in \eqref{proddz} in the form
\beql{diffo}
\si^{\tss\la}_{k\tss 0}(z)\ts\di_z^{\tss k}
+\si^{\tss\la}_{k1}(z)\ts\di_z^{\tss k-1}
+\dots+\si^{\tss\la}_{kk}(z)
\eeq
with
\ben
\si^{\tss\la}_{kl}(z)=\sum_{r\leqslant 0}\si^{\tss\la}_{kl}[r]z^{-r-l}.
\een
In particular,
\ben
\si^{\tss\la}_{11}[0]=\la_1+\dots+\la_{m+n},\qquad
\si^{\tss\la}_{11}[r]=e_{11}[r]+\dots+e_{m+n,m+n}[r],\qquad r<0,
\een
and the elements $\si^{\tss\la}_{kk}[r]$ satisfy
the recurrence relations
\ben
\si^{\tss\la}_{kk}[r]=\wt\si^{\ts\la}_{kk}[r]
+(\la_{m+n}+r+k-1)\ts\si^{\tss\la}_{k-1,k-1}[r]
+\sum_{p=r}^{-1}e_{m+n,m+n}[p]\tss\si^{\tss\la}_{k-1,k-1}[r-p]
\een
for $r<0$, and
\beql{recze}
\si^{\tss\la}_{kk}[0]=\wt\si^{\ts\la}_{kk}[0]
+(\la_{m+n}+k-1)\ts\si^{\tss\la}_{k-1,k-1}[0],
\eeq
where the elements $\wt\si^{\ts\la}_{kk}[r]$ are
defined by the expansion
\begin{multline}
\non
\sum_I \big({-}\di_z+\la_{i_1}z^{-1}+e_{i_1i_1}(z)_+\big)\dots
\big({-}\di_z+\la_{i_l}z^{-1}+e_{i_li_l}(z)_+\big)\\
{}\times\big(\ts\di_z+\la_{i_{\tss l+1}}z^{-1}
+e_{i_{\tss l+1}i_{\tss l+1}}(z)_+\big)
\dots
\big(\ts\di_z+\la_{i_k}z^{-1}+e_{i_ki_k}(z)_+\big)\\
{}=\wt\si^{\tss\la}_{k\tss 0}(z)\ts\di_z^{\tss k}
+\wt\si^{\tss\la}_{k1}(z)\ts\di_z^{\tss k-1}
+\dots+\wt\si^{\tss\la}_{kk}(z),
\end{multline}
summed over multisets
$I=\{m+n-1\geqslant
i_1\geqslant\dots\geqslant i_l\geqslant m+1>i_{l+1}>\dots>i_k\}$
with $l=0,\dots,k$
(we assume that $n\geqslant 1$ as the calculation in the case
$n=0$ is quite similar with some simplifications).

Observe that $\U(t^{-1}\h[t^{-1}])$ is a graded algebra
with the degree of the generator $e_{ii}[r]$
with $r<0$ equal to $-r$.
The above recurrence relations show that each
element $\si^{\tss\la}_{kk}[r]$ is homogeneous of degree $-r$.
Therefore, in order to complete the proof of the theorem
it will be sufficient to verify that each generator
$e_{ii}[r]$ can be expressed as a polynomial in the
elements $\si^{\tss\la}_{kk}[r]$. The algebraic independence
of these elements will then follow from the fact that
the homogeneous subspaces of $\U(t^{-1}\h[t^{-1}])$
are finite-dimensional.

We will use the reverse induction on $r$ and suppose that
$r\leqslant -1$. Using the above recurrence relations
we can write $\si^{\tss\la}_{kk}[r]$ as a homogeneous
polynomial of degree $-r$ in
the generators $e_{ii}[p]$ .
By the induction hypothesis, all generators $e_{ii}[p]$
with $p>r$ can be expressed as polynomials
in the $\si^{\tss\la}_{kk}[q]$. Furthermore, modulo
polynomials in those generators, each element
$\si^{\tss\la}_{kk}[r]$ with $k=1,\dots,m+n$
is given as a linear combination
of the generators $e_{ii}[r]$ with $i=1,\dots,m+n$.
Thus, we only need to show that the corresponding
matrix is nonsingular for any $r$.

As we only need to keep linear terms in the recurrence
relations, they are simplified to
\ben
\si^{\tss\la}_{kk}[r]=\wt\si^{\tss\la}_{kk}[r]
+(\la_{m+n}+r+k-1)\ts\si^{\tss\la}_{k-1,k-1}[r]
+e_{m+n,m+n}[r]\tss\si^{\tss\la}_{k-1,k-1}[0]
\een
for $r<0$. Using this together with \eqref{recze},
we derive that the desired linear combinations take
the form
\ben
\si^{\tss\la}_{kk}[r]=\sum_{i=1}^{m+n}\La_{ki}\ts e_{ii}[r],\qquad
k=1,\dots,m+n,
\een
where the coefficients $\La_{ki}$ are given by
\eqref{matele} and \eqref{matere}.
The determinant of the matrix $\La=[\La_{ki}]$
is found in Lemma~\ref{lem:det}. Recalling that $\la$
is a generic critical weight and applying \eqref{gener}
for the positive real roots of the form $\al=k\tss\de+\ve_i-\ve_j$
with $1\leqslant i<j\leqslant m$
we find that
$\la_i-\la_j+j-i+r\ne 0$ for all $r<0$.
Similarly, taking
$\al=k\tss\de+\ve_j-\ve_i$ with
$m+1\leqslant i<j\leqslant m+n$ and
$\al=k\tss\de+\ve_i-\ve_j$ with
$1\leqslant i\leqslant m<j\leqslant m+n$ we conclude
that for the respective values of $i$ and $j$
we have $\la_j-\la_i-j+i-r\ne 0$ for all $r<0$
and $\la_i+\la_j-i-j+2m+1\ne 0$. Hence,
the determinant of the matrix is nonzero.
This completes the proof.
\epf

Now consider the fields $h_{kl}(z)$ defined
in Corollary~\ref{cor:sop}.
Explicit formulas for these fields are found by using
the second relation in
Proposition~\ref{prop:extsh} which implies the identity
\ben
\bal
:\str^{}_{1,\dots,k}\ts H_k T_1(z)\dots T_k(z):
&=\sum_{I}\frac{1}{\al_{1}!\dots\al_{m}!}
\sum_{\si\in\Sym_k}\psi(\si,I,\si I)\\
{}&\times{}:\big(\di_z+\wh E(z)\big)_{i_1i_{\si(1)}}\dots
\big(\di_z+\wh E(z)\big)_{i_ki_{\si(k)}}:(-1)^{\ga(I,\si I)},
\eal
\een
where we keep
the notation of that proposition.
Introduce the
Fourier coefficients of the fields $h_{kl}(z)$ by the expansions
\ben
h_{kl}(z)=\sum_{r\in\ZZ}h_{kl}[r]z^{-r-l}.
\een

\bco\label{cor:gencom}
Suppose that $\la$ is a generic critical weight.
Then
\ben
M(\la)^{\wh\n}
=\CC[h_{kk}[r]\ |\  k=1,\dots,m+n,\  r<0\tss]\ts 1_{\la}
\een
so that any singular vector of $M(\la)$ can be expressed
as a unique polynomial in
the operators
$h_{kk}[r]$ with $k=1,\dots,m+n$ and $r<0$
applied to the highest vector.
\eco

\bpf
Let $h^{\la}_{kl}[r]$ denote the image of
$h_{kl}[r]\ts 1_{\la}$ under the isomorphism
\eqref{isomsv}.
As in the proof of Theorem~\ref{thm:gensv}, we will show
that the elements $h^{\la}_{kk}[r]$
with $k=1,\dots,m+n$ and $r<0$
are algebraically
independent generators of the algebra $\U(t^{-1}\h[t^{-1}])$.

Using generating series, we can
write relations \eqref{proddz} in the form
\begin{multline}
\non
\sum_{k=0}^{\infty}u^k\ts{}
:\str^{}_{1,\dots,k}\ts A_k T_1(z)\dots T_k(z):{}1_{\la}\mapsto\\
\Big(1+u\big(\di_z-\la_{m+n}z^{-1}
-e_{m+n,m+n}(z)_+\big)\Big)^{-1}\dots
\Big(1+u\big(\di_z-\la_{m+1}z^{-1}
-e_{m+1,m+1}(z)_+\big)\Big)^{-1}\\
{}\times{}\Big(1+u\big(\di_z+\la_{m}z^{-1}
+e_{m\tss m}(z)_+\big)\Big)\dots
\Big(1+u\big(\di_z+\la_{1}z^{-1}+e_{11}(z)_+\big)\Big).
\end{multline}
Hence, by Theorem~\ref{thm:mmmt},
\begin{multline}
\non
\sum_{k=0}^{\infty}u^k\ts{}
:\str^{}_{1,\dots,k}\ts H_k T_1(z)\dots T_k(z):{}1_{\la}\mapsto\\
\Big(1-u\big(\di_z+\la_{1}z^{-1}+e_{11}(z)_+\big)\Big)^{-1}\dots
\Big(1-u\big(\di_z+\la_{m}z^{-1}+e_{m\tss m}(z)_+\big)\Big)^{-1}\\
{}\times{}\Big(1+u\big({-}\di_z
+\la_{m+1}z^{-1}+e_{m+1,m+1}(z)_+\big)\Big)
\dots
\Big(1+u\big({-}\di_z+\la_{1}z^{-1}+e_{m+n,m+n}(z)_+\big)\Big).
\end{multline}
Therefore, taking the coefficient of $u^k$ we come to
the relation
\begin{multline}
\non
h^{\tss\la}_{k\tss 0}(z)\ts\di_z^{\tss k}
+h^{\tss\la}_{k1}(z)\ts\di_z^{\tss k-1}
+\dots+h^{\tss\la}_{kk}(z)\\
=\sum_I \big(\di_z+\la_{i_1}z^{-1}+e_{i_1i_1}(z)_+\big)\dots
\big(\di_z+\la_{i_l}z^{-1}+e_{i_li_l}(z)_+\big)\\
{}\times{}\big({-}\di_z+\la_{i_{\tss l+1}}z^{-1}
+e_{i_{\tss l+1}i_{\tss l+1}}(z)_+\big)
\dots
\big({-}\di_z+\la_{i_k}z^{-1}+e_{i_ki_k}(z)_+\big),
\end{multline}
summed over multisets
$I=\{i_1\leqslant\dots\leqslant i_l\leqslant m<i_{l+1}<\dots<i_k\}$
with $l=0,\dots,k$.
In particular, we have
\ben
h^{\tss\la}_{11}[0]=\la_1+\dots+\la_{m+n},\qquad
h^{\tss\la}_{11}[r]=e_{11}[r]+\dots+e_{m+n,m+n}[r],\qquad r<0,
\een
and the elements $h^{\tss\la}_{kk}[r]$ satisfy
the recurrence relations
\ben
h^{\tss\la}_{kk}[r]=\wt h^{\ts\la}_{kk}[r]
+(\la_{1}-r-k+1)\ts h^{\tss\la}_{k-1,k-1}[r]
+\sum_{p=r}^{-1}e_{11}[p]\tss h^{\tss\la}_{k-1,k-1}[r-p]
\een
for $r<0$, and
\ben
h^{\tss\la}_{kk}[0]=\wt h^{\ts\la}_{kk}[0]
+(\la_{1}-k+1)\ts h^{\tss\la}_{k-1,k-1}[0],
\een
where $\wt h^{\ts\la}_{kk}[r]$ are the elements
of the algebra $\U(t^{-1}\h_{m-1|n}[t^{-1}])$
defined as the coefficients
in the expansion
of the corresponding differential operator
associated with the natural subalgebra
$\h_{m-1|n}\subset \h$ generated by
the elements $e_{ii}[r]$ with $2\leqslant i\leqslant m+n$.
Arguing as in the proof of Theorem~\ref{thm:gensv}, we
will keep only linear terms in the expansion of
$h^{\tss\la}_{kk}[r]$ so that the recurrence
relations simplify to
\ben
h^{\tss\la}_{kk}[r]=\wt h^{\ts\la}_{kk}[r]
+(\la_{1}-r-k+1)\ts h^{\tss\la}_{k-1,k-1}[r]
+e_{11}[r]\tss h^{\tss\la}_{k-1,k-1}[0]
\een
for $r<0$. Observe now that if we change the notation by
setting $\la_i:=-\la_{m+n-i+1}$,
$e_{ii}[r]:=e_{m+n-i+1,m+n-i+1}[r]$ for $r<0$, and
\ben
\si^{\tss\la}_{kk}[r]:=\begin{cases}(-1)^{k+1}\ts h^{\tss\la}_{kk}[r]
\qquad&\text{if}\quad r<0,\\
(-1)^{k}\ts h^{\tss\la}_{kk}[r]
\qquad&\text{if}\quad r=0,
\end{cases}
\een
then the recurrence relations would take exactly the same form
as those in the proof of Theorem~\ref{thm:gensv}
for the elements $\si^{\tss\la'}_{kk}[r]$
associated with the Lie superalgebra $\wh\gl_{n|m}$
and the highest weight
\ben
\la'=(-\la_{m+n},\dots,-\la_{m+1},-\la_m,\dots,-\la_1).
\een
It remains to verify that
the corresponding determinant of Lemma~\ref{lem:det}
does not vanish under the specialization at $\la'$.
This value is easily found and it turns out to be given
by the same product as in Lemma~\ref{lem:det}, with
each linear factor replaced by its negative.
It was checked in the proof of Theorem~\ref{thm:gensv}
that each factor is nonzero.
\epf

Applying the Newton identity \eqref{newtonthm},
we can also calculate the images
under the isomorphism $\HC$ of
the singular vectors of the Verma module $M(\la)$
associated with the normally ordered
traces of powers of the matrix $T(z)$;
see Corollary~\ref{cor:sop}. Using the notation
introduced in the proofs of Theorem~\ref{thm:gensv}
and Corollary~\ref{cor:gencom} set
\ben
h^{\tss\la}_k(z)=h^{\tss\la}_{k\tss 0}(z)\ts\di_z^{\tss k}
+h^{\tss\la}_{k1}(z)\ts\di_z^{\tss k-1}
+\dots+h^{\tss\la}_{kk}(z)
\een
and
\ben
\si^{\tss\la}_k(z)=\si^{\tss\la}_{k\tss 0}(z)\ts\di_z^{\tss k}
+\si^{\tss\la}_{k1}(z)\ts\di_z^{\tss k-1}
+\dots+\si^{\tss\la}_{kk}(z).
\een
Denote by $s^{\tss\la}_{kl}(z)$ the image of the series
$s_{kl}(z)\ts 1_{\la}$ under the isomorphism $\HC$
defined in \eqref{isomsv}, so that
\ben
:\str\ts T(z)^k:{} 1_{\la}=s^{\tss\la}_{k\tss 0}(z)\ts
\di_z^{\tss k}+s^{\tss\la}_{k1}(z)\ts\di_z^{\tss k-1} +\dots
+s^{\tss\la}_{kk}(z).
\een

\bco\label{cor:eipos}
The series $s^{\tss\la}_{kl}(z)$ can be found from
the relations
\ben
s^{\tss\la}_{k\tss 0}(z)\ts
\di_z^{\tss k}+s^{\tss\la}_{k1}(z)\ts\di_z^{\tss k-1} +\dots
+s^{\tss\la}_{kk}(z)=
\sum_{l=0}^{k-1}(-1)^l\tss(l+1)\ts h^{\tss\la}_{k-l-1}(z)\ts
\si^{\tss\la}_{l+1}(z).
\een
\eco

\bpf
This is immediate from Theorem~\ref{thm:ident}.
\epf

\end{document}